\documentclass[12pt,a4paper]{article}
\usepackage[a4paper]{geometry}
\usepackage{amssymb,latexsym,amsmath,amsfonts,amsthm}
\usepackage{graphicx}
\usepackage{epstopdf}
\usepackage{tikz}
\usepackage{pgflibraryshapes}
\usetikzlibrary{arrows,decorations.markings}
\usepackage{overpic}
\usepackage{fullpage}
\usepackage{comment,verbatim}
\usepackage{hyperref}
\usepackage[title,titletoc]{appendix}
\usepackage{float}
\usepackage{appendix}

\usepackage{caption}
\usepackage{subfigure}
\usepackage{authblk}
\usepackage{arydshln}
\usepackage{enumerate}

\usepackage[thicklines]{cancel}

\newtheorem{thm}{Theorem}[section]

\newtheorem{pro}[thm]{Proposition}

\theoremstyle{remark}
\newtheorem{rem}[thm]{Remark}

\numberwithin{equation}{section}

\def\Im {\mathop{\rm Im}\nolimits}
\def\arg {\mathop{\rm arg}\nolimits}
\def\Re {\mathop{\rm Re}\nolimits}

\begin{document}

\title{On the asymptotics of real solutions for the Painlev\'{e} I equation}

\author{Wen-Gao Long$^1$,\ \ \ \ \  Jun Xia$^{2,}$\footnote{Corresponding author}}
\date{}
\maketitle
\begin{center}
{$^1$\footnotesize
{\textit{School of Mathematics and Computational Science, Hunan University of Science and Technology,\\
Xiangtan, 411201, China,}
\\
\texttt{longwg@hnust.edu.cn}

\medskip

$^2$\footnotesize
\textit{School of Mathematics and Systems Science, Guangdong Polytechnic Normal University, \\
Guangzhou 510665, China,}
\\
\texttt{junxiacqu@foxmail.com}}
\\
}
\end{center}

\medskip
\medskip

\begin{abstract}
In this paper, we revisit the asymptotic formulas of real Painlev\'e I transcendents as the independent variable tends to negative infinity, which were initially derived by Kapaev with the complex WKB method. Using the Riemann-Hilbert method, we improve the error estimates of the oscillatory type asymptotics and provide precise error estimates of the singular type asymptotics. We also establish the corresponding asymptotics for the associated Hamiltonians of real Painlev\'e I transcendents. In addition, two typos in the mentioned asymptotic behaviors in literature are corrected. \\
\newline
  \textbf{2010 mathematics subject classification:} 33E17; 34A30; 34E05; 34M55; 41A60
\newline
  \textbf{Keywords and phrases:} The Painlev\'{e} I equation, asymptotic expansions, Riemann-Hilbert approach.
\end{abstract}

\maketitle
\numberwithin{equation}{section}


\section{Introduction}
Painlev\'{e} equations are a group of six nonlinear second-order ordinary differential equations
that possess the Painlev\'{e} property. This property requires that the movable singularities of the solutions must be poles, not branch points or essential singularities.
In general, the solutions of these equations do not have explicit expressions. Consequently,
asymptotic analysis is a common tool to study the behavior of solutions for Painlev\'{e} equations.

In this paper, we are concerned with the Painlev\'e I equation
\begin{equation}\label{PI}
y''(x)=6y^2(x)+x,
\end{equation}
and we concentrate on the asymptotic behaviors as  $x\to-\infty$ of real solutions to \eqref{PI}.
Using the method of dominant balance (cf. \cite[Chapter 3]{BO}), it is easy to see that there are two classes of behaviors
\begin{align}
y(x)\sim \left(-\frac{x}{6}\right)^{\frac12}, \qquad or \qquad
y(x)\sim -\left(-\frac{x}{6}\right)^{\frac12},
\end{align}
as $x\to-\infty$. For the second case, one can further show that as $x\to-\infty$,
\begin{align}\label{y}
y(x)=-\left(-\frac{x}{6}\right)^{\frac12}+\frac{d}{(-x)^{\frac18}}
\cos\left(\frac{4\cdot24^{\frac14}}5(-x)^{\frac54}-\frac{5d^2}8\ln(-x)+\chi\right)
+O(x^{-\frac58}),
\end{align}
where $d$ and $\chi$ are constants; see \cite{JK,QL}. It should be noted that the asymptotic formula \eqref{y} also appears in the NIST handbook of mathematical functions (see \cite[Eq. 32.11.1]{NIST}).

Indeed, Holmes and Spence \cite{HS} proved that there exist exactly three types
of real solutions by studying
a boundary value problem for \eqref{PI}. Subsequently, with the help of the isomonodromy approach \cite{FIKN}, Kapaev \cite{Kap1} obtained the following classification of asymptotic behaviors as $x\to-\infty$ of real Painlev\'e I transcendents, in terms of the Stokes multiplier $s_2$ (see Section \ref{Section2} for definition).
\begin{enumerate}[(A)]
  \item When $|s_2|<1$, we have
  \begin{align}\label{y1}
y(x)\sim-\left(-\frac{x}{6}\right)^{\frac12}+\frac{d}{(-24x)^{\frac18}}
\cos\left(\frac{4\cdot24^{\frac14}}5(-x)^{\frac54}-\frac{5d^2}8\ln(-x)+\chi\right),
\end{align}
where
\begin{equation}\label{con-for-1}
\left\{\begin{aligned}
d^2&=-\frac1\pi\ln(1-|s_2|^2),\\
\chi&=\arg s_2-\left(\frac{19}{8}\ln2+\frac58\ln3\right)d^2-\arg\Gamma\left(-\frac{id^2}{2}\right)
+\frac{3\pi}4.
\end{aligned}\right.
\end{equation}
  \item When $|s_2|=1$, we have
  \begin{align}\label{y2}
y(x)=\widehat{y}(x)+\frac{s_1-s_{-1}}{4\cdot24^{\frac14}\sqrt{\pi}}(-x)^{-\frac18}
e^{-\frac{4\cdot24^{\frac14}}5(-x)^{\frac54}}\left[1+O(x^{-\frac54})\right],
\end{align}
where $\widehat{y}(x)=(-\frac{x}{6})^{\frac12}\left[1+O(x^{-\frac52})\right]$ as $x\to-\infty$ is the solution to \eqref{PI} with the Stokes multipliers $s_2=s_{-2}=i,~s_{-1}=s_1=\frac{i}{2}$.
  \item When $|s_2|>1$, we have
\begin{align}\label{y3}
y(x)\sim-\left(-\frac{x}{6}\right)^{\frac12}+\frac{(-x)^{\frac12}}
{\frac{\sqrt{6}}{3}\sin^2\left(\frac{2\cdot24^{\frac14}}5(-x)^{\frac54}
+\frac{5q}8\ln(-x)+\rho\right)},
\end{align}
where
\begin{equation}\label{con-for-2}
\left\{\begin{aligned}
\rho&=\frac1{2\pi}\ln(|s_2|^2-1),\\
\sigma&=\frac12\arg s_2+\left(\frac{19}{8}\ln2+\frac58\ln3\right)\rho
+\frac12\arg\Gamma\left(\frac{1}{2}-i\rho\right)-\frac\pi4.
\end{aligned}\right.
\end{equation}
\end{enumerate}

The type (A), (B) and (C) solutions are called oscillatory, separatrix and singular solutions, respectively. In \cite{Kap2}, with the Riemann-Hilbert approach, Kapaev found all separatrix (also known as \emph{tronqu\'ee}) solutions of \eqref{PI}. Recently, Dea\~{n}o \cite{Deano} extends Kapaev's results by calculating explicitly exponentially small corrections with respect to the tronqu\'ee solutions. However, to the best of our knowledge, the error estimate in asymptotic formula \eqref{y3} and the asymptotics for the Hamiltonians in case (A) and (C) have not been provided in the existing literature.

In this paper, by performing a Deift-Zhou steepest descent analysis \cite{Deift,DZ1,DZ2} to the Riemann-Hilbert problem of \eqref{PI}, we obtain the following asymptotics for the Painlev\'e I transcendent $y(x)$ and its Hamiltonian $\mathcal{H}(x)$ (see \eqref{Hamilton} for definition).

\begin{thm}\label{thm}
If $|s_2|<1$, we have, as $x\to-\infty$,
\begin{align}\label{y-asymp1}
y(x)&=-\left(-\frac{x}{6}\right)^{\frac12}+\frac{\sqrt{a}}{(-24x)^{\frac18}}
\cos\left(\frac{4\cdot24^{\frac14}}5(-x)^{\frac54}-\frac{5a}8\ln(-x)+\phi\right)
+O(x^{-\frac34}),\\
\mathcal{H}(x)&=-4\left(-\frac{x}{6}\right)^{\frac32}
+a\left(-\frac{3x}2\right)^{\frac14}+\frac{\sqrt{a}}{(-24x)^{\frac38}}\nonumber
\sin\left(\frac{4\cdot24^{\frac14}}5(-x)^{\frac54}-\frac{5a}8\ln(-x)+\phi\right)\\
&\qquad+O(x^{-1}),  \label{H-asymp1}
\end{align}
where
\begin{equation}
\left\{\begin{aligned}
a&=-\frac1\pi\ln(1-|s_2|^2),\\
\phi&=\arg s_2-\left(\frac{19}{8}\ln2+\frac58\ln3\right)a-\arg\Gamma\left(-\frac{ia}{2}\right)
-\frac\pi4.
\end{aligned}\right.
\end{equation}
If $|s_2|>1$, we have, as $x\to-\infty$,
\begin{align}\label{y-asymp2}
y(x)&=-\left(-\frac{x}{6}\right)^{\frac12}+\frac{(-x)^{\frac12}}
{\frac{\sqrt{6}}{3}\sin^2\left(\frac{2\cdot24^{\frac14}}5(-x)^{\frac54}
+\frac{5b}8\ln(-x)+\psi\right)}+O(x^{-\frac34}),\\
\mathcal{H}(x)&=-4\left(-\frac{x}{6}\right)^{\frac32}
-b\left(-24x\right)^{\frac14}-\left(-\frac{3x}2\right)^{\frac14}\nonumber
\cot\left(\frac{2\cdot24^{\frac14}}5(-x)^{\frac54}+\frac{5b}8\ln(-x)+\psi\right)\\
&\qquad+O(x^{-1}),  \label{H-asymp2}
\end{align}
where the error terms hold uniformly for $x$ bounded away from the singularities appearing in the denominator and
\begin{equation}
\left\{\begin{aligned}
b&=\frac1{2\pi}\ln(|s_2|^2-1),\\
\psi&=\frac12\arg s_2+\left(\frac{19}{8}\ln2+\frac58\ln3\right)b
+\frac12\arg\Gamma\left(\frac{1}{2}-ib\right)+\frac\pi4.
\end{aligned}\right.
\end{equation}
\end{thm}

\begin{figure}[h]
\centering
\includegraphics[scale=0.57]{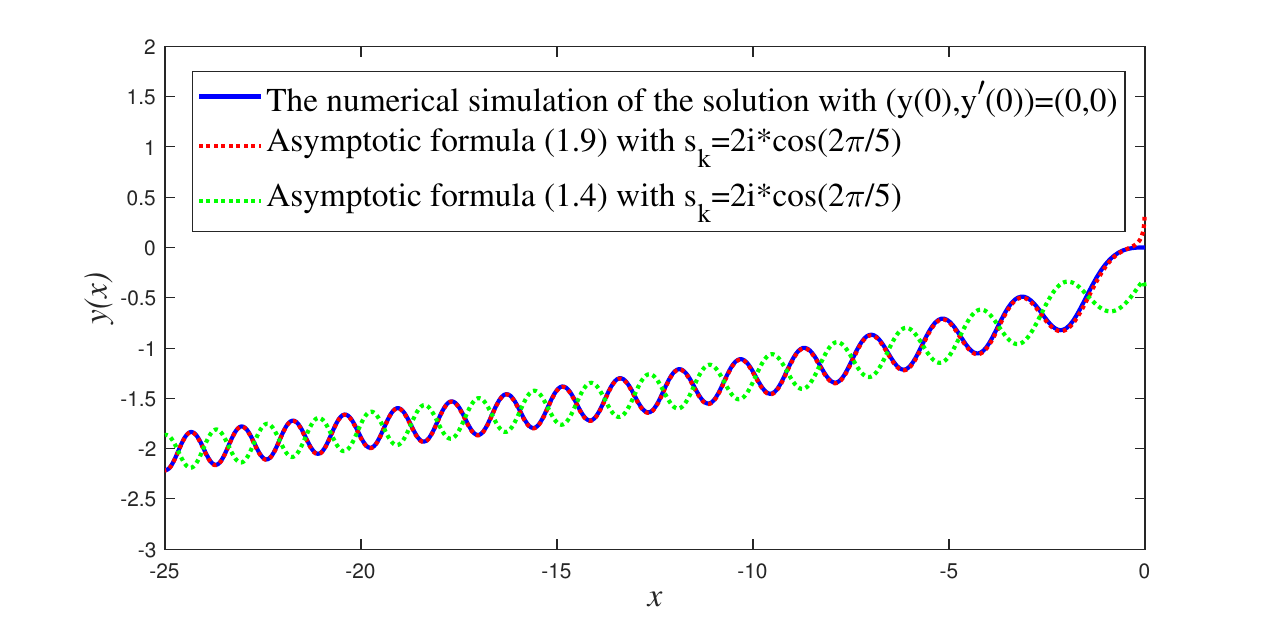}
\caption{A comparison of the numerical and asymptotical solutions for the Painlev\'{e} I equation with $(y(0),y'(0))=(0,0)$.}\label{Osci-solu}
\end{figure}

\begin{figure}[h]
\centering
\includegraphics[scale=0.59]{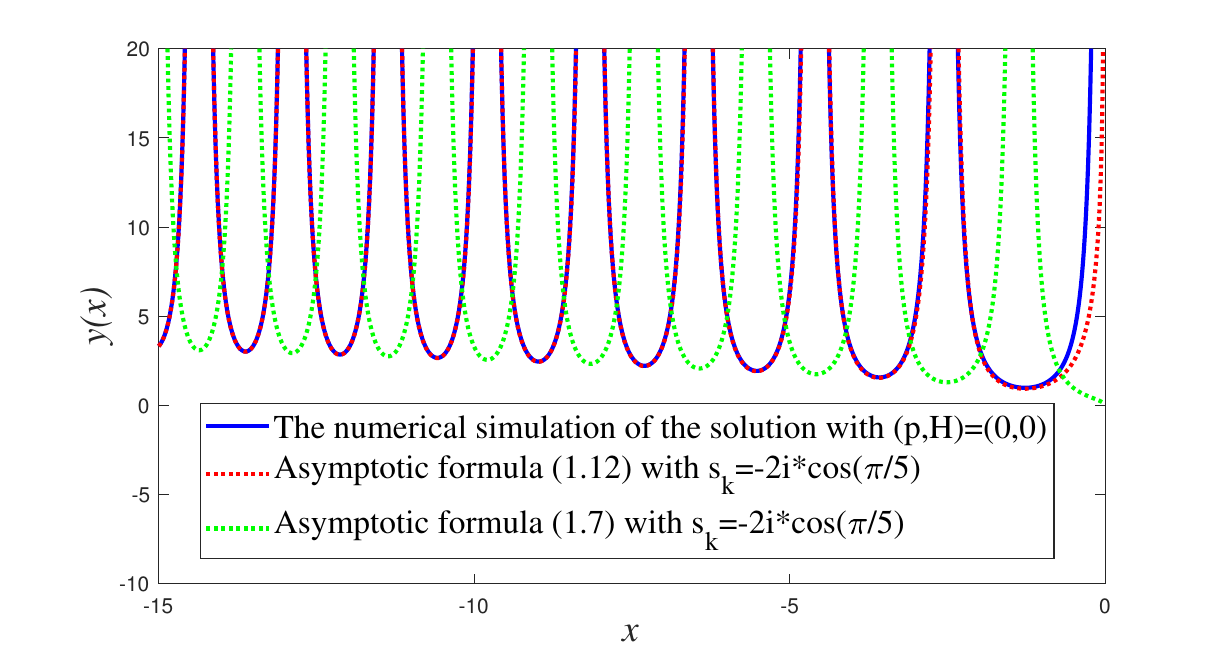}
\caption{A comparison of the numerical and asymptotical solutions for the Painlev\'{e} I equation with $(p,H)=(0,0)$.}\label{singular-solu}
\end{figure}

\begin{rem}
In the asymptotic formula \eqref{y-asymp1}, we improve the error estimate $O(x^{-5/8})$ in \eqref{y} to $O(x^{-3/4})$. In the asymptotic formula \eqref{y-asymp2}, we provide a rigorous error estimate $O(x^{-3/4})$. Moreover, we find that the angle $\frac{3\pi}{4}$ in \eqref{con-for-1} and  $-\frac\pi4$ in \eqref{con-for-2} should be corrected by $-\frac\pi4$ and $\frac\pi4$, respectively.
The first author (W.-G. Long) of this paper apologies for committing the same typos in the Introduction of his other papers \cite{LongLi, LongLLZ, LongLiWang}, even though these typos do not affect the correctness of the main results of these papers. In Figures \ref{Osci-solu} and \ref{singular-solu}, we have verified \eqref{y-asymp1} and \eqref{y-asymp2} numerically in special cases.
\end{rem}

\begin{rem}
One should note from \cite{Kitaev} that the Stokes multipliers of the special Painlev\'{e} I transcendent with $(y(0),y'(0))=(0,0)$ are $s_{k}=2i\cos\left(\frac{2\pi}{5}\right), k=0, \pm 1,\pm 2$.
In Appendix \ref{appendix2}, we calculate the Stokes multipliers of the special solution with $(p,H)=(0,0)$ by using the complex WKB method and find that $s_{k}=-2i\cos\left(\frac{\pi}{5}\right), k=0, \pm 1,\pm 2$, where $p$ is the pole location of the solution and $H$ is a free parameter in the Laurent series
\begin{equation}
y(x)=\frac{1}{(x-p)^{2}}-\frac{p}{10}(x-p)^2-\frac{1}{6}(x-p)^3+H(x-p)^4+\cdots.
\end{equation}
With this result, we successfully verify our asymptotic formula \eqref{y-asymp2} numerically for this special solution.
\end{rem}

\begin{figure}[h]
\centering
\includegraphics[scale=0.50]{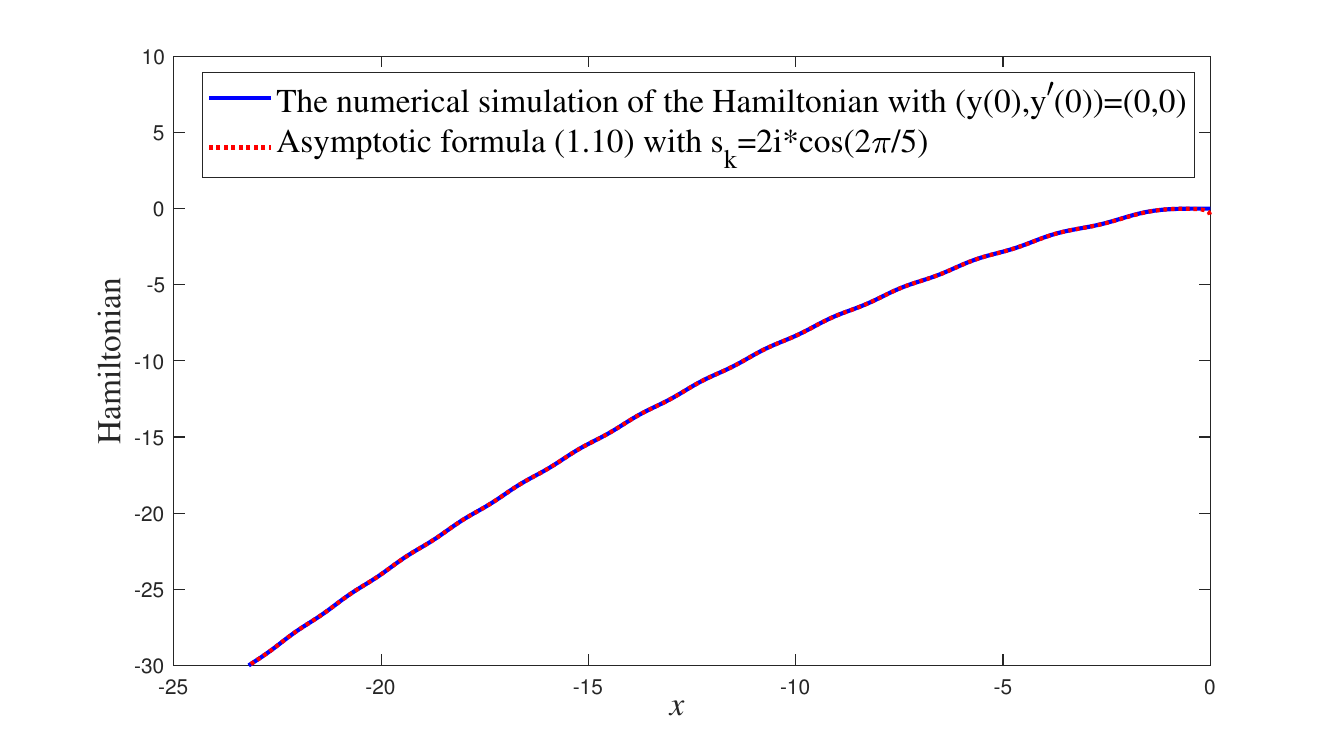}
\caption{A comparison of the numerical and asymptotical solutions for the Hamiltonian with $(y(0),y'(0))=(0,0)$.}\label{oscillate-H}
\end{figure}

\begin{figure}[h]
\centering
\includegraphics[scale=0.51]{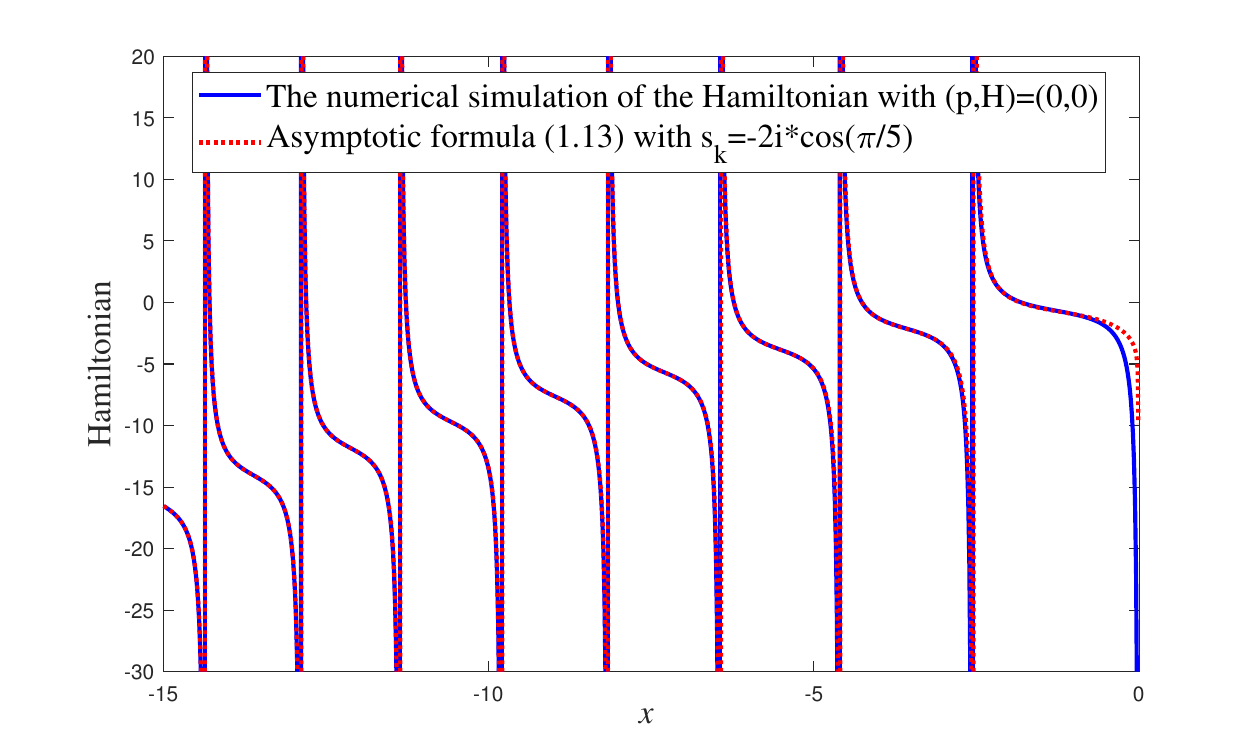}
\caption{A comparison of the numerical and asymptotical solutions for the Hamiltonian with $(p,H)=(0,0)$.}\label{singular-H}
\end{figure}

\begin{rem}
In Figure \ref{oscillate-H} and Figure \ref{singular-H}, we numerically verify the asymptotic behaviors of the Hamiltonian $\mathcal{H}(x)$ in \eqref{H-asymp1} and \eqref{H-asymp2}, respectively.
\end{rem}

The rest of the paper is organized as follows. In Section \ref{Section2}, we formulate the Riemann-Hilbert problem of the Painlev\'e I equation \eqref{PI}. Due to the different case: (1) $|s_2|<1$, (2) $|s_2|>1$, we apply the Deift-Zhou approach to analyze the Riemann-Hilbert problem for $\Psi(z,x)$ as $x\to-\infty$ in Section \ref{Section3} and \ref{Section4}, respectively. In the final section \ref{Section5}, we prove Theorem \ref{thm}. For the convenience of the reader, we put the parabolic cylinder and Airy parametrix model used in the asymptotic analysis, and the computation of Stokes multipliers of a special solution to Painlev\'e I in the Appendix.

\section{Riemann-Hilbert problem}\label{Section2}
\begin{figure}[t]
  \centering
  \includegraphics[width=10cm,height=8.8cm]{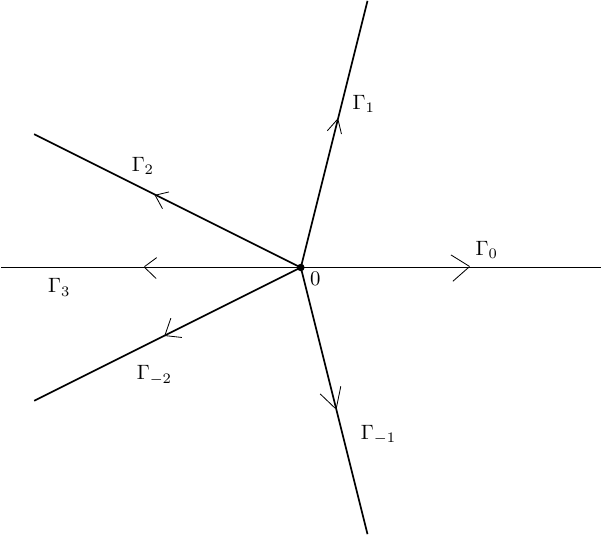}\\
  \caption{The jump contour $\Gamma$.}\label{PIjump}
\end{figure}

In this section, we describe the Riemann-Hilbert (RH) problem of the Painlev\'e I equation \eqref{PI}, which was first introduced in \cite{Kap2}.
We seek a $2\times2$ matrix-valued function $\Psi(z,x)$ satisfying the following properties.
\begin{enumerate}[(a)]
\item $\Psi(z,x)$ is analytic for $z\in \mathbb{C}\setminus\Gamma$, where $\Gamma=\cup_{k=-2}^{2}\Gamma_{k}$ with
\begin{align*}
\Gamma_0=\mathbb{R}_+,\quad \Gamma_3=e^{\pi i}\mathbb{R}_+,\quad \Gamma_{\pm1}=e^{\pm\frac{2\pi i}{5}}\mathbb{R}_+,\quad \Gamma_{\pm2}=e^{\pm\frac{4\pi i}{5}}\mathbb{R}_+,
\end{align*}
see Figure \ref{PIjump} for an illustration.
\item $\Psi(z,x)$ satisfies the jump relation
\begin{align}\label{Psi-jump}
\Psi_{+}(z,x)=\Psi_{-}(z,x)S_{k},\ \ \ \ \ z\in\Gamma_{k},
\end{align}
   where the Stokes matrices $\{S_k\}$ take the forms
   \begin{equation}
  \begin{aligned}
  &S_{0}=\begin{pmatrix}1 & 0\\s_{0} & 1\end{pmatrix},\ \ \ \ \ \
  S_{1}= \begin{pmatrix}1 & s_{1}\\0 & 1\end{pmatrix},\ \ \ \ \ \ \ \
  S_{2}=\begin{pmatrix}1 & 0\\s_{2} & 1\end{pmatrix},\\
  &S_{3}=\begin{pmatrix}0 & -i\\-i & 0\end{pmatrix},\ \ \ \
  S_{-1}= \begin{pmatrix}1 & s_{-1}\\0 & 1\end{pmatrix},\ \ \ \
  S_{-2}=\begin{pmatrix}1 & 0\\s_{-2} & 1\end{pmatrix}.
  \end{aligned}
  \end{equation}
  Here, the Stokes multipliers $\{s_k\}$ satisfy the restricted conditions
\begin{align}\label{Stokes-eq1}
1+s_{k}s_{k+1}=-is_{k+3},\qquad s_{k+5}=s_{k},\qquad k\in\mathbb{Z}.
\end{align}
\item $\Psi(z;x)$ satisfies the following asymptotic behavior as $z\to\infty$
  \begin{align}\label{Psi-at-infty}
  \Psi(z,x)=z^{\frac{\sigma_{3}}{4}}\frac{1}{\sqrt{2}}
  \begin{pmatrix}1 & 1\\ 1 & -1\end{pmatrix}\left[I+\sum_{k=1}^{\infty}\Psi_k(x)z^{-\frac{k}{2}}\right]
  e^{(\frac{4}{5}z^{\frac{5}{2}}+xz^{\frac{1}{2}})\sigma_{3}},
 \end{align}
  where the branches of powers of $z$ are chosen such that $\arg z\in(-\pi,\pi)$.
\end{enumerate}

It should be mentioned that in the RH analysis throughout, we will use frequently the matrix notations
$$
\sigma_1=\begin{pmatrix} 0 & 1 \\ 1 & 0\end{pmatrix},\qquad
\sigma_2=\begin{pmatrix} 0 & -i \\ i & 0\end{pmatrix},\qquad
\sigma_3=\begin{pmatrix} 1 & 0 \\ 0 & -1\end{pmatrix}.
$$

From \cite[Equations (2.13)-(2.15)]{Kap2}, we obtain the following expressions for the Painlev\'e I transcendent $y(x)$ and the associated Hamiltonian
\begin{align}\label{Hamilton}
\mathcal{H}(x)=\frac12[y'(x)]^2-2y(x)^3-xy(x).
\end{align}
\begin{pro}
Let $\Psi_k:=\Psi_k(x)$, $k=1,2$ be the coefficient in the expansion \eqref{Psi-at-infty}. We have
  \begin{align}\label{Psi1Psi2}
 \Psi_1=\begin{pmatrix} -\mathcal{H}(x) &0\\0&\mathcal{H}(x)  \end{pmatrix},\qquad
 \Psi_2=\frac12\begin{pmatrix} \mathcal{H}(x)^2 &y(x)\\y(x)&\mathcal{H}(x)^2\end{pmatrix}.
  \end{align}
Hence, we have
\begin{align}\label{y-express}
y(x)&=2\Psi_{2,12}=2\Psi_{2,21},\\
\mathcal{H}(x)&=-\Psi_{1,11}=\Psi_{1,22},\label{H-express}
\end{align}
where $\Psi_{k,ij}$ denotes the $(i,j)$ entry of $\Psi_k$.
\end{pro}

By \eqref{Stokes-eq1}, it follows that
\begin{align}\label{Stokes-eq2}
1+s_2s_{-2}=-is_0,\qquad s_1=\frac{i-s_{-2}}{1+s_2s_{-2}},\qquad s_{-1}=\frac{i-s_{2}}{1+s_2s_{-2}},
\end{align}
if $1+s_2s_{-2}\neq0$. For the real solutions of the Painlev\'e I equation, we further have (see \cite{Kap1})
\begin{align}\label{Stokes-eq3}
\bar{s}_2=-s_{-2},
\end{align}
which implies from \eqref{Stokes-eq2} that
\begin{align}\label{Stokes-eq4}
\bar{s}_0=-s_{0},\qquad \bar{s}_1=-s_{-1},
\end{align}
and
\begin{align}\label{Stokes-eq5}
-is_0=1+s_2s_{-2}=1-|s_2|^2=1-|s_{-2}|^2\in\mathbb{R}.
\end{align}

\section{Asymptotics of $\Psi$ for $|s_2|<1$}\label{Section3}
In this section, we first perform the Deift-Zhou nonlinear steepest descent analysis to the RH problem for $\Psi(z,x)$ with $0<|s_2|<1$. Then, the reduced case $s_2=0$ will be considered in the end of this section.

\subsection{Scaling transformation}
Assume in what follows that $x<0$. We define
\begin{align}\label{Psi-to-A}
A(z)=|x|^{-\frac{\sigma_3}{8}}\Psi(|x|^{\frac12}z,x).
\end{align}
Then, $A(z)$ satisfies the following RH problem.
\begin{enumerate}[(a)]
\item $A(z)$ is analytic for $z\in \mathbb{C}\setminus\Gamma$.

\item $A(z)$ satisfies the same jump relation as $\Psi(z,x)$; see \eqref{Psi-jump}.

\item As $z\to\infty$, we have
  \begin{align}\label{A-at-infty}
  A(z)=z^{\frac{\sigma_{3}}{4}}\frac{1}{\sqrt{2}}
  \begin{pmatrix}1 & 1\\ 1 & -1\end{pmatrix}\left[I+\sum_{k=1}^{\infty}
  \Psi_k(x)(|x|^{\frac12}z)^{-\frac{k}{2}}\right]
  e^{|x|^{\frac{5}{4}}\theta(z)\sigma_{3}},
 \end{align}
  where
  \begin{align}\label{theta}
\theta(z)=\frac{4}{5}z^{\frac{5}{2}}-z^{\frac{1}{2}}.
  \end{align}
\end{enumerate}

\subsection{The $g$-function transformation}
\begin{figure}[t]
  \centering
  \includegraphics[width=9cm,height=7cm]{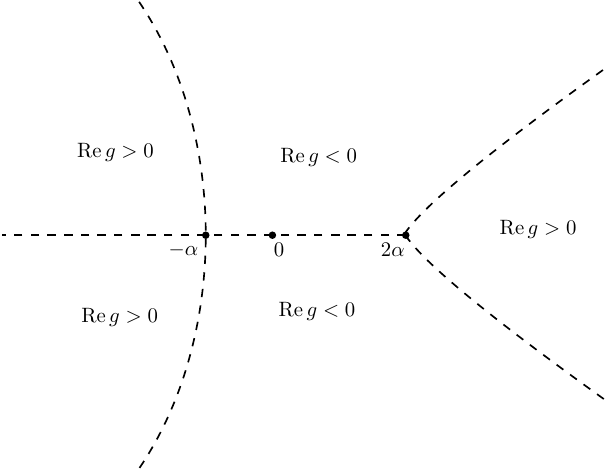}\\
  \caption{The sign of $\Re g(z)$ in the complex plane. The dashed curves stand for $\Re g(z)=0$.}\label{g-sign}
\end{figure}

To normalize the behavior of $A(z)$ at infinity, we introduce the $g$-function
\begin{align}\label{g}
g(z)=\frac{4}{5}\left(z-2\alpha\right)^{\frac32}
\left(z+3\alpha\right),\qquad \alpha:=\frac{\sqrt{6}}{6}.
\end{align}
where the branch of $(z-2\alpha)^{\frac32}$ is taken such that $\arg(z-2\alpha)\in(-\pi,\pi)$. It is clear that
\begin{align}\label{g-expansion}
g(z)=\theta(z)+\frac{2\alpha}{3}z^{-\frac12}+O(z^{-\frac32}),\qquad z\to\infty,
\end{align}
and $g(z)$ has saddle point at $z=-\alpha$. The sign of $\Re g(z)$ in the complex plane is depicted in Figure \ref{g-sign}.

The next transformation is defined by
\begin{align}\label{A-to-B}
B(z)=\begin{pmatrix}1&\frac{2\alpha}{3}|x|^{\frac54}-\Psi_{1,11}|x|^{-\frac14}\\
0&1\end{pmatrix}
A(z)e^{-|x|^{\frac{5}{4}}g(z)\sigma_3}.
\end{align}
We come to the following RH problem for $B(z)$.
\begin{enumerate}[(a)]
\item $B(z)$ is analytic for $z\in \mathbb{C}\setminus\Gamma$.

\item $B(z)$ satisfies the jump relation
\begin{align}\label{B-jump}
B_+(z)=B_-(z)\left\{\begin{aligned}
&e^{|x|^{\frac{5}{4}}g(z)\sigma_3}S_ke^{-|x|^{\frac{5}{4}}g(z)\sigma_3}, &&z\in\Gamma_{k},~k=\pm1,\pm2,\\
&-i\sigma_1, &&z\in\Gamma_3,\\
&e^{|x|^{\frac{5}{4}}g_-(z)\sigma_3}S_0e^{-|x|^{\frac{5}{4}}g_+(z)\sigma_3}, &&z\in(0,2\alpha),\\
&e^{|x|^{\frac{5}{4}}g(z)\sigma_3}S_0e^{-|x|^{\frac{5}{4}}g(z)\sigma_3}, &&z\in(2\alpha,+\infty).
\end{aligned}\right.
\end{align}

\item As $z\to\infty$, we have
  \begin{align}\label{B-at-infty}
  B(z)=\left[I+\frac{B_1}{z}+O(z^{-2})\right]z^{\frac{\sigma_{3}}{4}}\frac{1}{\sqrt{2}}
  \begin{pmatrix}1 & 1\\ 1 & -1\end{pmatrix},
 \end{align}
 where
 \begin{align}\label{B1-21}
&B_{1,21}=\Psi_{1,11}|x|^{-\frac14}-\frac{2\alpha}{3}|x|^{\frac54},\\ &B_{1,22}=\left(\Psi_{2,11}-\Psi_{2,12}\right)|x|^{-\frac12}
-\frac{2\alpha}{3}\Psi_{1,11}|x|+\frac{2\alpha^2}{9}|x|^{\frac52}.\label{B1-22}
 \end{align}
\end{enumerate}

\subsection{Contour deformation}\label{contour-defor}
\begin{figure}[t]
  \centering
  \includegraphics[width=10cm,height=7cm]{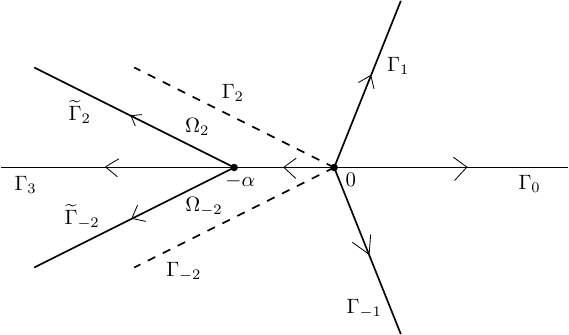}\\
  \caption{The contour $\Gamma_{C}$.}\label{deform-1}
\end{figure}

In view of Figure \ref{g-sign}, we see from \eqref{B-jump} that the jumps on $\Gamma_{\pm2}$ do not tend to the identity matrix as $x\to-\infty$. To overcome this, we define the transformation
\begin{equation}\label{B-to-C}
C(z)=B(z)\left\{\begin{aligned}
&e^{|x|^{\frac{5}{4}}g(z)\sigma_3}S_2^{-1}e^{-|x|^{\frac{5}{4}}g(z)\sigma_3}, && z\in \Omega_2,\\
&e^{|x|^{\frac{5}{4}}g(z)\sigma_3}S_{-2}e^{-|x|^{\frac{5}{4}}g(z)\sigma_3}, && z\in \Omega_{-2},\\
&I, && \mathrm{elsewhere}.
\end{aligned}\right.
\end{equation}
Thus, we have the following RH problem for $C(z)$.
\begin{enumerate}[(a)]
\item $C(z)$ is analytic for $z\in \mathbb{C}\setminus\Gamma_C$, where $\Gamma_C$ is shown in Figure \ref{deform-1}.

\item $C(z)$ satisfies the jump relation
\begin{align}\label{C-jump}
C_+(z)=C_-(z)\left\{\begin{aligned}
&e^{|x|^{\frac{5}{4}}g(z)\sigma_3}S_ke^{-|x|^{\frac{5}{4}}g(z)\sigma_3}, &&z\in\Gamma_{k},~k=\pm1,\\
&e^{|x|^{\frac{5}{4}}g(z)\sigma_3}S_ke^{-|x|^{\frac{5}{4}}g(z)\sigma_3}, &&z\in\widetilde\Gamma_{k},~k=\pm2,\\
&-i\sigma_1, &&z\in(-\infty,-\alpha),\\
&e^{|x|^{\frac{5}{4}}g_-(z)\sigma_3}S_2(-i\sigma_1)S_{-2}
e^{-|x|^{\frac{5}{4}}g_+(z)\sigma_3}, &&z\in(-\alpha,0),\\
&e^{|x|^{\frac{5}{4}}g_-(z)\sigma_3}S_0e^{-|x|^{\frac{5}{4}}g_+(z)\sigma_3}, &&z\in(0,2\alpha),\\
&e^{|x|^{\frac{5}{4}}g(z)\sigma_3}S_0e^{-|x|^{\frac{5}{4}}g(z)\sigma_3}, &&z\in(2\alpha,+\infty).
\end{aligned}\right.
\end{align}

\item As $z\to\infty$, we have
\begin{align}\label{C-at-infty}
C(z)=\left[I+\frac{B_1}{z}+O(z^{-2})\right]z^{\frac{\sigma_{3}}{4}}\frac{1}{\sqrt{2}}
\begin{pmatrix}1 & 1\\ 1 & -1\end{pmatrix},
 \end{align}
 where $B_1$ is given in \eqref{B-at-infty}.
\end{enumerate}

\begin{figure}[t]
  \centering
  \includegraphics[width=11cm,height=6.8cm]{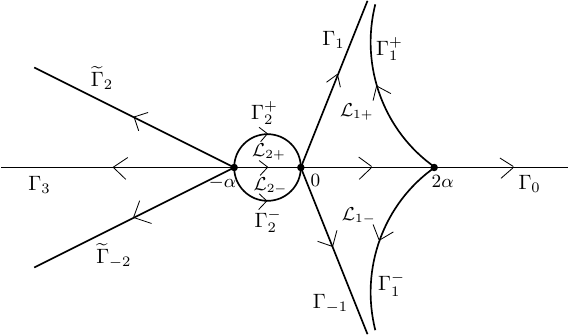}\\
  \caption{The contour $\Gamma_{D}$.}\label{deform-2}
\end{figure}

On the other hand, since $g_{\pm}(z)$ is purely imaginary on $(-\infty,\alpha)$, the diagonal entries of jumps on $(-\alpha,2\alpha)$ are highly oscillating for large $|x|$. To turn the oscillations into exponential decays, we open lens around the segment $(-\alpha,2\alpha)$ according to the sign of $\Re g(z)$ depicted in Figure \ref{g-sign}. This process depends on the following matrix decomposition:
\begin{align}\label{decompo-1}
&[S_2(-i\sigma_1)S_{-2}]^{-1}=\begin{pmatrix} 1&-is_2s_0^{-1}\\ 0&1 \end{pmatrix}
\begin{pmatrix}0&-s_0^{-1}\\ s_0&0 \end{pmatrix}
\begin{pmatrix} 1&-is_{-2}s_0^{-1}\\ 0&1 \end{pmatrix},\\
&S_0=\begin{pmatrix} 1&s_0^{-1}\\ 0&1 \end{pmatrix}
\begin{pmatrix}0&-s_0^{-1}\\ s_0&0 \end{pmatrix}
\begin{pmatrix} 1&s_0^{-1}\\ 0&1 \end{pmatrix}.\label{decompo-2}
\end{align}
Owing to \eqref{decompo-1} and \eqref{decompo-2}, we define
\begin{equation}\label{C-to-D}
D(z)=C(z)\left\{\begin{aligned}
&\begin{pmatrix} 1&-s_0^{-1}e^{2|x|^{\frac{5}{4}}g(z)}\\ 0&1 \end{pmatrix}, &&z\in\mathcal{L}_{1+},\\
&\begin{pmatrix} 1&s_0^{-1}e^{2|x|^{\frac{5}{4}}g(z)}\\ 0&1 \end{pmatrix}, &&z\in\mathcal{L}_{1-},\\
&\begin{pmatrix} 1&is_{-2}s_0^{-1}e^{2|x|^{\frac{5}{4}}g(z)}\\ 0&1 \end{pmatrix}, &&z\in\mathcal{L}_{2+},\\
&\begin{pmatrix} 1&-is_2s_0^{-1}e^{2|x|^{\frac{5}{4}}g(z)}\\ 0&1 \end{pmatrix}, &&z\in\mathcal{L}_{2-},
\end{aligned}\right.
\end{equation}
where the regions $\mathcal{L}_{1\pm}$, $\mathcal{L}_{2\pm}$ are depicted in Figure \ref{deform-2}. Therefore, we arrive at the following RH problem for $D(z)$.
\begin{enumerate}[(a)]
\item $D(z)$ is analytic for $z\in \mathbb{C}\setminus\Gamma_D$, where $\Gamma_D$ is shown in Figure \ref{deform-2}.

\item $D(z)$ satisfies the jump relations
\begin{align}\label{D-jump1}
D_+(z)=D_-(z)\left\{\begin{aligned}
&\begin{pmatrix} 1&-is_{\mp2}s_0^{-1}e^{2|x|^{\frac{5}{4}}g(z)}\\ 0&1 \end{pmatrix}, &&z\in\Gamma_{\pm1},\\
&\begin{pmatrix} 1&0\\ s_{\pm2}e^{-2|x|^{\frac{5}{4}}g(z)}&1 \end{pmatrix}, &&z\in\widetilde\Gamma_{\pm2},\\
&\begin{pmatrix} 1&-s_0^{-1}e^{2|x|^{\frac{5}{4}}g(z)}\\ 0&1 \end{pmatrix}, &&z\in\Gamma^{+}_1\cup\Gamma^-_1,\\
&\begin{pmatrix} 1&-is_{\mp2}s_0^{-1}e^{2|x|^{\frac{5}{4}}g(z)}\\ 0&1 \end{pmatrix}, &&z\in\Gamma^{\pm}_2,
\end{aligned}\right.
\end{align}
and
\begin{align}\label{D-jump2}
D_+(z)=D_-(z)\left\{\begin{aligned}
&\begin{pmatrix}0&-i\\ -i&0 \end{pmatrix}, &&z\in(-\infty,-\alpha),\\
&\begin{pmatrix}0&-s_0^{-1}\\ s_0&0 \end{pmatrix}, &&z\in(-\alpha,2\alpha),\\
&\begin{pmatrix}1&0\\ s_0e^{-2|x|^{\frac{5}{4}}g(z)}&1 \end{pmatrix}, &&z\in(2\alpha,+\infty).
\end{aligned}\right.
\end{align}

\item As $z\to\infty$, we have
  \begin{align}\label{D-at-infty}
D(z)=\left[I+\frac{B_1}{z}+O(z^{-2})\right]z^{\frac{\sigma_{3}}{4}}\frac{1}{\sqrt{2}}
\begin{pmatrix}1 & 1\\ 1 & -1\end{pmatrix},
 \end{align}
 where $B_1$ is given in \eqref{B-at-infty}.
\end{enumerate}

Again, it is readily seen from the signs of $\Re g(z)$ shown in Figure \ref{g-sign} that as $x\to-\infty$, the jumps of $D(z)$ tend to identity matrix except on $(-\infty,2\alpha)$.

\subsection{Global parametrix}
In this section, we construct the global parametrix. By \eqref{D-jump2} and \eqref{D-at-infty}, we consider the following model RH problem for a function $N(z)$.
\begin{enumerate}[(a)]
\item $N(z)$ is analytic for $z\in \mathbb{C}\setminus(-\infty,2\alpha)$, where $(-\infty,2\alpha)$ is oriented from left to right.

\item $N(z)$ satisfies the jump relation
\begin{align}\label{N-jump}
N_+(z)=N_-(z)\left\{\begin{aligned}
&\begin{pmatrix}0&i\\ i&0 \end{pmatrix}, &&z\in(-\infty,-\alpha),\\
&\begin{pmatrix}0&-s_0^{-1}\\ s_0&0 \end{pmatrix}, &&z\in(-\alpha,2\alpha).
\end{aligned}\right.
\end{align}

\item As $z\to\infty$, we have
  \begin{align}\label{N-at-infty}
N(z)=\left[I+O(z^{-1})\right]z^{\frac{\sigma_{3}}{4}}\frac{1}{\sqrt{2}}
  \begin{pmatrix}1 & 1\\ 1 & -1\end{pmatrix}.
 \end{align}
\end{enumerate}

It is straightforward to verify that
\begin{align}\label{N}
N(z)=\begin{pmatrix}1&2\sqrt{3\alpha}\,i\nu\\0&1 \end{pmatrix}
(z-2\alpha)^{\frac{\sigma_{3}}4}
\frac{1}{\sqrt{2}}\begin{pmatrix}1 & 1\\ 1 & -1\end{pmatrix}d(z)^{\sigma_3},
\end{align}
where $d(z)$ is the Szeg\H{o} function
\begin{align}\label{d}
d(z)&=\exp\left(-\frac{\ln(-is_0)}{2\pi}
\sqrt{z-2\alpha}\int^{2\alpha}_{-\alpha}\frac{d\xi}
{\sqrt{2\alpha-\xi}\,(\xi-z)}\right)\nonumber\\
&=\left(\frac{i\sqrt{z-2\alpha}+\sqrt{3\alpha}}
{i\sqrt{z-2\alpha}-\sqrt{3\alpha}}\right)^{\nu}
\end{align}
with
\begin{align}\label{nu}
\nu=-\frac1{2\pi i}\ln(-is_0)=-\frac1{2\pi i}\ln(1-|s_2|^2)\in i\mathbb{R}.
\end{align}
The branches of the multi-valued functions in \eqref{N} and \eqref{d} are fixed by the conditions
$$
\arg(z-2\alpha)\in(-\pi,\pi), \qquad \arg\left(\frac{i\sqrt{z-2\alpha}+\sqrt{3\alpha}}
{i\sqrt{z-2\alpha}-\sqrt{3\alpha}}\right)\in(-\pi,\pi).
$$
A direct calculation gives
\begin{align}\label{N-expansion}
N(z)=\left[I+\frac{N_1}{z}+O(z^{-2})\right]
z^{\frac{\sigma_{3}}4}\frac{1}{\sqrt{2}}\begin{pmatrix}1 & 1\\ 1 & -1\end{pmatrix},\qquad z\to\infty,
\end{align}
where
\begin{align}\label{N1}
N_1=\begin{pmatrix}
6\alpha\nu^2-\frac{\alpha}{2} & 2\sqrt{3}\alpha^{\frac32}i\nu(1-4\nu^2)\\
-2\sqrt{3\alpha}\,i\nu & \frac{\alpha}{2}-6\alpha\nu^2
\end{pmatrix}.
\end{align}

\subsection{Local parametrices}
In this section, we construct relevant local parametrices in the neighborhood $U(z_0)=\{z\in\mathbb{C}:|z-z_0|<\delta\}$ of the special point $z_0\in\{-\alpha,0,2\alpha\}$.

Firstly, we focus on the saddle point $z_0=-\alpha$. In view of the RH problem for $D(z)$, we consider the following local RH problem.
\begin{enumerate}[(a)]
\item $P^{(-)}(z)$ is analytic for $z\in U(-\alpha)\setminus\Gamma_D$.

\item $P^{(-)}(z)$ satisfies the same jump relation as $D(z)$ on $\Gamma_D\cap U(-\alpha)$.

\item For $z\in\partial U(-\alpha)$, we have
  \begin{align}\label{matching-}
  P^{(-)}(z)N(z)^{-1}=I+O(|x|^{-\frac{5}{8}}),\qquad x\to-\infty.
 \end{align}
\end{enumerate}

Let us define a conformal mapping
\begin{align}\label{lambda-1}
\lambda_1(z)=2\left\{\begin{aligned}
&\sqrt{g_+(-\alpha)-g(z)}, \quad &&\Im z>0,\\
&\sqrt{g_+(-\alpha)+g(z)}, \quad &&\Im z<0,
\end{aligned}\right.
\end{align}
where the branch of the square roots are chosen by the condition
\begin{align}\label{lambda-1-asymp}
\lambda_1(z)=e^{-\frac{\pi i}{4}}2(3\alpha)^{\frac14}(z+\alpha)\left[1-\frac{z+\alpha}{18\alpha}
+O\left((z+\alpha)^2\right)\right],\qquad z\to-\alpha,
\end{align}
and
\begin{align}\label{g(-alpha)}
g_+(-\alpha)=-\frac{4}{5}\sqrt{3\alpha}\,i.
\end{align}
Then, we construct $P^{(-)}(z)$ as
\begin{align}\label{P-}
P^{(-)}(z)=E^{(-)}(z)\Phi^{(\mathrm{PC})}\left(|x|^{\frac58}\lambda_1(z)\right)
s_2^{-\frac{\sigma_3}{2}}h_1^{\frac{\sigma_3}{2}}M(z)e^{-|x|^{\frac54}g(z)\sigma_3},
\end{align}
where $\Phi^{(\mathrm{PC})}$ is the parabolic cylinder parametrix model with the parameter $\nu$ given in \eqref{nu} (see Appendix \ref{PCP}), $h_{1}$ is defined in \eqref{h0h1} and
\begin{align}\label{E-}
E^{(-)}(z)&=N(z)M(z)^{-1}h_1^{-\frac{\sigma_3}{2}}s_2^{\frac{\sigma_3}{2}}
e^{-|x|^{\frac54}g_+(-\alpha)\sigma_3}
\left(|x|^{\frac58}\lambda_1(z)\right)^{\nu\sigma_3}2^{-\frac{\sigma_3}{2}}
\begin{pmatrix}\lambda_1(z) &1 \\ 1 & 0\end{pmatrix}
\end{align}
with
\begin{align}\label{M}
M(z)=\left\{\begin{aligned}
&i\sigma_1,&&\Im z>0,\\
&I,&&\Im z<0.\end{aligned}\right.
\end{align}
With \eqref{Stokes-eq2} and \eqref{PC-jump}, we check directly that $P^{(-)}(z)$ fulfills the same jump relation as $D(z)$ on $\Gamma_D\cap U(-\alpha)$ and $E^{(-)}(z)$ is analytic in $U(-\alpha)$. In addition, we obtain from \eqref{P-}, \eqref{E-} and \eqref{PC-at-infty} that as $x\to-\infty$,
\begin{align}\label{PN-1}
P^{(-)}(z)N(z)^{-1}&=N(z)M(z)^{-1}\left[I+\frac{L(z)}{|x|^{\frac58}}
+O(|x|^{-\frac54})\right]M(z)N(z)^{-1}\nonumber\\
&=I+O(|x|^{-\frac58}),
\end{align}
where
\begin{align}\label{L}
L(z)=\begin{pmatrix} 0 & \frac{\nu s_2|x|^{\frac54\nu}\lambda_1(z)^{2\nu-1}}
{h_1e^{2|x|^{\frac54}g_+(-\alpha)}}\\ \frac{h_1e^{2|x|^{\frac54}g_+(-\alpha)}}
{s_2|x|^{\frac54\nu}\lambda_1(z)^{2\nu+1}} & 0
\end{pmatrix}.
\end{align}
Notice that we have made use of the fact that $\nu$ and $g_+(-\alpha)$ are purely imaginary (see \eqref{nu} and \eqref{g(-alpha)}).

For the local parametrix near the origin $z_0=0$, we consider the following local RH problem for a function $P^{(0)}(z)$.
\begin{enumerate}[(a)]
\item $P^{(0)}(z)$ is analytic for $z\in U(0)\setminus\Gamma_D$.

\item $P^{(0)}(z)$ satisfies the same jump relation as $D(z)$ on $\Gamma_D\cap U(0)$.

\item For $z\in\partial U(0)$, we have
  \begin{align}\label{matching0}
  P^{(0)}(z)N(z)^{-1}=I+O(|x|^{-\infty}),\qquad x\to-\infty.
 \end{align}
\end{enumerate}

The solution to above problem is elementary. First, we note from the real solution condition \eqref{Stokes-eq3} that
\begin{align}\label{s2s-2}
s_2=re^{i\beta},\qquad s_{-2}=-re^{-i\beta},
\end{align}
where $r:=|s_2|$, $\beta:=\arg s_2$. Define
\begin{equation}\label{F}
F(\lambda)=e^{\lambda\sigma_3}\left\{\begin{aligned}
&I,&&\arg\lambda\in\textstyle(\frac{\pi}{4},\frac{3\pi}{4})
\cup(-\frac{3\pi}{4},-\frac{\pi}{4}),\\
&\begin{pmatrix}1& ir\\0&1  \end{pmatrix}, &&\arg\lambda\in\textstyle(\frac{3\pi}{4},\frac{5\pi}{4}),\\
&\begin{pmatrix}1& 0\\-ir&1  \end{pmatrix}, &&\arg\lambda\in\textstyle(-\frac{\pi}{4},\frac{\pi}{4}).
\end{aligned}\right.
\end{equation}
Then, we see that $F(\lambda)$ solves the following RH problem.
\begin{enumerate}[(a)]
\item $F(\lambda)$ is analytic for $\lambda\in \mathbb{C}\setminus\{\lambda:\arg\lambda=\pm\frac\pi4,\pm\frac{3\pi}{4}\}$, where rays are oriented from $0$ to $\infty$.

\item $F(\lambda)$ satisfies the jump relation
\begin{align}\label{F-jump}
F_+(\lambda)=F_-(\lambda)\left\{\begin{aligned}
&\begin{pmatrix}1&0\\ ir&1 \end{pmatrix}, &&\arg\lambda=\frac\pi4,\\
&\begin{pmatrix}1&ir\\ 0&1 \end{pmatrix}, &&\arg\lambda=\frac{3\pi}4,\\
&\begin{pmatrix}1&-ir\\ 0&1 \end{pmatrix}, &&\arg\lambda=-\frac{3\pi}4,\\
&\begin{pmatrix}1&0\\ -ir&0 \end{pmatrix}, &&\arg\lambda=-\frac\pi4.
\end{aligned}\right.
\end{align}

\item As $\lambda\to\infty$, we have
  \begin{align}\label{F-at-infty}
F(\lambda)=\left[I+O(\lambda^{-\infty})\right]e^{\lambda\sigma_3}.
 \end{align}
\end{enumerate}

Introduce the conformal mapping
\begin{align}\label{lambda-2}
\lambda_2(z)&=g_-(0)\pm g(z),\qquad \pm\Im z>0,\nonumber\\
&=e^{\frac{\pi i}{2}}2\sqrt{2\alpha}z(1+o(1)),\qquad z\to0.
\end{align}
We construct $P^{(0)}(z)$ as
\begin{align}\label{P0}
P^{(0)}(z)=E^{(0)}(z)F\left(|x|^{\frac54}\lambda_2(z)\right)
e^{\frac{i\beta\sigma_3}{2}}Q(z)
s_0^{\frac{\sigma_3}{2}}e^{-|x|^{\frac54}g(z)\sigma_3},
\end{align}
where
\begin{align}\label{E0}
E^{(0)}(z)=N(z)s_0^{-\frac{\sigma_3}{2}}Q(z)^{-1}e^{-\frac{i\beta\sigma_3}{2}}
e^{-|x|^{\frac54}g_-(0)\sigma_3}
\end{align}
and
\begin{align}\label{Q}
Q(z)=\left\{\begin{aligned}
&I, &&\Im z>0,\\
&\begin{pmatrix}0&1\\-1&0 \end{pmatrix}, &&\Im z<0.
\end{aligned}\right.
\end{align}
Using \eqref{N-jump}, \eqref{s2s-2} and \eqref{F-jump}, it is direct to see from \eqref{P0}-\eqref{Q} that $P^{(0)}(z)$ satisfies the same jump relation as $D(z)$ on $\Gamma_D\cap U(0)$ and $E^{(0)}(z)$ is analytic in $U(0)$. Moreover, it follows from \eqref{F-at-infty} and the fact that $g_-(0)$ is purely imaginary that the matching condition \eqref{matching0} is also satisfied.

Finally, we concentrate on the branch point $z_0=2\alpha$ and we consider the following local RH problem for some function $P^{(+)}(z)$.
\begin{enumerate}[(a)]
\item $P^{(+)}(z)$ is analytic for $z\in U(2\alpha)\setminus\Gamma_D$.

\item $P^{(+)}(z)$ satisfies the same jump relation as $D(z)$ on $\Gamma_D\cap U(2\alpha)$.

\item For $z\in\partial U(2\alpha)$, we have
  \begin{align}\label{matching+}
  P^{(+)}(z)N(z)^{-1}=I+O(|x|^{-\frac{5}{4}}),\qquad x\to-\infty.
 \end{align}
\end{enumerate}

As usual, in the neighborhood $U(2\alpha)$, we define a conformal mapping
\begin{align}\label{lambda-3}
\lambda_3(z)=\left(\frac32 g(z)\right)^{\frac23}=(6\alpha)^{\frac23}(z-2\alpha)(1+o(1)),\qquad z\to2\alpha.
\end{align}
Then the desired parametrix near $z_0=2\alpha$ is defined by
\begin{align}\label{P+}
P^{(+)}(z)=E^{(+)}(z)\Phi^{(\mathrm{Ai})}\left(|x|^{\frac56}\lambda_3(z)\right)\sigma_1
s_0^{\frac{\sigma_3}{2}}e^{-|x|^{\frac54}g(z)\sigma_3},
\end{align}
where $\Phi^{(\mathrm{Ai})}$ is the Airy parametrix model (see Appendix \ref{AP}) and
\begin{align}\label{E+}
E^{(+)}(z)=N(z)s_0^{-\frac{\sigma_3}{2}}\sigma_1\frac{1}{\sqrt{2}} \begin{pmatrix}1 & -i\\ -i &1\end{pmatrix}\left(|x|^{\frac56}\lambda_3(z)\right)^{\frac{\sigma_{3}}{4}}.
\end{align}

\subsection{Small norm problem}
\begin{figure}[t]
  \centering
  \includegraphics[width=11cm,height=6.8cm]{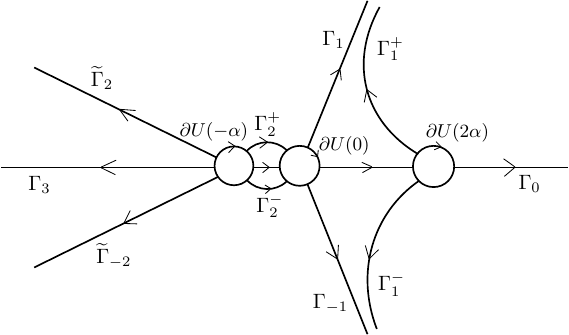}\\
  \caption{The contour $\Gamma_{R}$.}\label{Gamma-R}
\end{figure}

Using the model functions $N(z)$ given in \eqref{N}, $P^{(-)}(z)$ given in \eqref{P-}, $P^{(0)}(z)$ given in \eqref{P0} and $P^{(+)}(z)$ given in \eqref{P+}, we define the final transformation as
\begin{align}\label{D-to-R}
R(z)=D(z)\left\{\begin{aligned}
&P^{(-)}(z)^{-1}, &&z\in U(-\alpha),\\
&P^{(0)}(z)^{-1}, &&z\in U(0),\\
&P^{(+)}(z)^{-1}, &&z\in U(2\alpha),\\
&N(z)^{-1}, &&\mathrm{elsewhere}.
\end{aligned}\right.
\end{align}
As a result, we have the following RH problem for $R(z)$.
\begin{enumerate}[(a)]
\item $R(z)$ is analytic for $z\in \mathbb{C}\setminus\Gamma_R$, where $\Gamma_R$ is illustrated in Figure \ref{Gamma-R}.

\item $R(z)$ satisfies the jump relation $R_+(z)=R_-(z)J_R(z)$
\begin{equation}\label{R-jump}
J_R(z)=\left\{\begin{aligned}
&P^{(-)}(z)N(z)^{-1}, &&z\in\partial U(-\alpha),\\
&P^{(0)}(z)N(z)^{-1}, &&z\in\partial U(0),\\
&P^{(+)}(z)N(z)^{-1}, &&z\in\partial U(2\alpha),\\
&N(z)J_D(z)N(z)^{-1}, &&z\in\Gamma_R\setminus(\partial U(-\alpha)\cup\partial U(0)\cup\partial U(2\alpha)),
\end{aligned}\right.
\end{equation}
where $J_D(z)$ denotes the relevant matrix in \eqref{D-jump1} and \eqref{D-jump2}.

\item As $z\to\infty$, we have
  \begin{align}\label{R-at-infty}
 R(z)=I+\frac{R_1}{z}+O(z^{-2}).
 \end{align}
\end{enumerate}

From the matching conditions \eqref{matching-}, \eqref{matching0} and \eqref{matching+}, we see from \eqref{R-jump} that the jump matrix $J_R(z)$ satisfy the estimation
\begin{align}\label{Rjump-estimate}
J_R(z)=I+O(|x|^{-\frac58}),\qquad x\to-\infty.
\end{align}
With the standard small norm theory \cite{Deift}, we have
\begin{align}\label{R-estimate}
R(z)=I+O(|x|^{-\frac58}),\qquad x\to-\infty,
\end{align}
uniformly for $z\in\mathbb{C}\setminus\Gamma_R$.

\subsection{The reduced case $|s_2|=0$}\label{sec:szero}
In the reduced case $|s_2|=0$, we obtain from \eqref{Stokes-eq5} that $s_0=i$. One sees from the RH problem for $D(z)$ in Section \ref{contour-defor} that there are no jumps on $\widetilde{\Gamma}_{\pm2}\cup\Gamma_2^{\pm}\cup\Gamma_{\pm1}$. The subsequent asymptotic analysis is almost same as the case $0<|s_2|<1$.
First of all, we solve a RH problem with jump matrix $i\sigma_1$ on $(-\infty,2\alpha)$. The solution to this RH problem can be constructed as \eqref{N} by taking $\nu=0$ therein. Next, we need to build a local parametrix near $z=2\alpha$. The desired parametrix is defined as \eqref{P+} by taking $|s_2|=0$ therein. Similar to \eqref{D-to-R}, by defining $R(z)$ as the ratio of the RH problem for $D(z)$ and the parametrices, it is readily seen that $R(z)$ satisfies the estimation $R(z)=I+O(|x|^{-5/4})$ as $x\to-\infty$.

\section{Asymptotics of $\Psi$ for $|s_2|>1$}\label{Section4}
In this section, we will perform the Deift-Zhou nonlinear steepest descent analysis to the RH problem for $\Psi(z,x)$ with $|s_2|>1$. As in the case $0<\rho<1$, it includes a series of invertible transformations, and the first four transformations $\Psi(z,x)\to{A(z)}\to{B(z)}\to{C(z)}\to{D(z)}$ are exactly the same as  those defined in  \eqref{Psi-to-A}, \eqref{A-to-B}, \eqref{B-to-C} and\eqref{C-to-D}, respectively. The main differences lie in the constructions of the global and local parametrices. Indeed, by ignoring the exponentially small terms of $J_D(z)$, we are led to the same RH problem for $N(z)$. Since $|s_2|>1$, a solution is given by \eqref{N} but with $\nu$ therein replaced by
\begin{equation}\label{nu0}
\nu=-\frac{1}{2}-\frac{1}{2\pi{i}}\ln(|s_2|^2-1)
:=-\frac{1}{2}+\nu_0.
\end{equation}
It is obvious that $\nu_0\in{i}\mathbb{R}$. If we continue to construct the local parametrix $P^{(-)}$ as in \eqref{P-} and \eqref{E-}, it comes out that $P^{(-)}(z)N(z)^{-1}$ does not tend to the identity matrix as $t\to+\infty$ for $z\in \partial U(-\alpha)$; see \eqref{PN-1} with $\nu$ given by \eqref{nu0}. This implies that the matching condition \eqref{matching-} is no longer valid. To overcome this, our idea is to make some modifications of the global and local parametrices.

\subsection{Global parametrix}
We modify the global parametrix by imposing specific singularity at $z_0=-\alpha$. More precisely, we consider the following RH problem.
\begin{enumerate}[(a)]
\item $\widehat{N}(z)$ is analytic for $z\in \mathbb{C}\setminus(-\infty,2\alpha)$, where $(-\infty,2\alpha)$ is oriented from left to right.

\item $\widehat{N}(z)$ satisfies the same jump relation as $N(z)$; see \eqref{N-jump}.

\item As $z\to\infty$, we have
  \begin{align}\label{Nhat-at-infty}
\widehat{N}(z)=\left[I+O(z^{-1})\right]z^{\frac{\sigma_{3}}{4}}\frac{1}{\sqrt{2}}
  \begin{pmatrix}1 & 1\\ 1 & -1\end{pmatrix}.
 \end{align}

\item As $z\to-\alpha$, $\widehat{N}(z)$ behaves like
\begin{align}
\widehat{N}(z)=O((z+\alpha)^{-\frac32}).
\end{align}
\end{enumerate}

We look for a solution in the form
\begin{align}\label{Nhat}
\widehat{N}(z)=\left(\frac{X}{z+\alpha}+I\right)N(z),\qquad \det\left(\frac{X}{z+\alpha}+I\right)=1,
\end{align}
where $N(z)$ is given via \eqref{N} and \eqref{d}. The $z$-independent matrix $X$ is to be determined. It is easy to see that
\begin{align}\label{Nhat-expansion}
\widehat{N}(z)=\left[I+\frac{X+N_1}{z}+O(z^{-2})\right]
z^{\frac{\sigma_{3}}{4}}\frac{1}{\sqrt{2}}
  \begin{pmatrix}1 & 1\\ 1 & -1\end{pmatrix},\qquad z\to\infty,
\end{align}
where $N_1$ is given in \eqref{N1}.

\subsection{Local parametrices}
Near the saddle point $z_0=-\alpha$, we instead consider the following local RH problem.
\begin{enumerate}[(a)]
\item $\widehat{P}^{(-)}(z)$ is analytic for $z\in U(-\alpha)\setminus\Gamma_D$.

\item $\widehat{P}^{(-)}(z)$ satisfies the same jump relation as $D(z)$ on $\Gamma_D\cap U(-\alpha)$.

\item For $z\in\partial U(-\alpha)$, we have
  \begin{align}\label{matching-m}
  \widehat{P}^{(-)}(z)\widehat{N}(z)^{-1}=I+O(|x|^{-\frac{5}{4}}),\qquad x\to-\infty.
 \end{align}
\end{enumerate}

We define the local parametrix $\widehat{P}^{(-)}(z)$ as
\begin{align}\label{P-m}
\widehat{P}^{(-)}(z)=\widehat{E}^{(-)}(z)
\Phi^{(\mathrm{PC})}\left(|x|^{\frac58}\lambda_1(z)\right)
s_2^{-\frac{\sigma_3}{2}}h_1^{\frac{\sigma_3}{2}}M(z)e^{-|x|^{\frac54}g(z)\sigma_3},
\end{align}
where $\Phi^{(\mathrm{PC})}$ is the parabolic cylinder parametrix model with the parameter $\nu$ given in \eqref{nu0} (see Appendix \ref{PCP}), $\lambda_1(z)$ is the conformal mapping \eqref{lambda-1} and
\begin{align}\label{E-m}
\widehat{E}^{(-)}(z)=\widehat{N}(z)M(z)^{-1}h_1^{-\frac{\sigma_3}{2}}
&s_2^{\frac{\sigma_3}{2}}e^{-|x|^{\frac54}g_+(-\alpha)\sigma_3}
\left(|x|^{\frac58}\lambda_1(z)\right)^{\nu\sigma_3}\nonumber\\
&\times\begin{pmatrix}
1 & 0\\ -\frac{1}{|x|^{\frac58}\lambda_1(z)} & 1
\end{pmatrix}2^{-\frac{\sigma_3}{2}}
\begin{pmatrix}\lambda_1(z) &1 \\ 1 & 0\end{pmatrix}.
\end{align}

Using the jumps \eqref{D-jump1}, \eqref{D-jump2} and \eqref{PC-jump}, we see that $\widehat{P}^{(-)}(z)$ has the same jump matrices as $D(z)$ on $\Gamma_D\cap U(-\alpha)$ and $\widehat{E}^{(-)}(z)$ is analytic in $U(-\alpha)\setminus\{-\alpha\}$. In addition, with the asymptotic behavior \eqref{PC-at-infty}, it is easily seen that the matching condition \eqref{matching-m} is fulfilled. To make sure that $\widehat{E}^{(-)}(z)$ is also analytic at the isolated point $z=-\alpha$, we have to make an appropriate choice of $X$ in \eqref{Nhat}. Indeed, by calculating the Laurent expansion at $z=-\alpha$ of $\widehat{E}^{(-)}(z)$ using \eqref{N}, \eqref{d}, \eqref{lambda-1-asymp}, \eqref{P-m} and \eqref{E-m}, we have
\begin{align}\label{Ehat}
\widehat{E}^{(-)}(z)=\left(\frac{X}{z+\alpha}+I\right)
\left(\frac{H}{z+\alpha}+I\right)Y(z),
\end{align}
where $Y(z)$ is analytic near $z_0=-\alpha$ and
\begin{align}\label{H}
H=\frac{6\alpha c}{1-c}\begin{pmatrix}2\nu_0 &-4\sqrt{3\alpha}i\nu_0^2\\
-i/\sqrt{3\alpha} & -2\nu_0\end{pmatrix},
\end{align}
with
\begin{align}\label{c}
c=c(x)=\frac{h_1e^{2|x|^{\frac54}g_+(-\alpha)}}{s_2|x|^{\frac54\nu_0}e^{\frac{3\pi i\nu_0}{2}}2^{6\nu_0}(3\alpha)^{\frac{5\nu_0}{2}}}.
\end{align}
Since $H^2=0$, it follows from \eqref{Ehat} that the $X$ is determined by
\begin{align}\label{X}
X=-H=\frac{6\alpha c}{1-c}\begin{pmatrix}-2\nu_0 & 4\sqrt{3\alpha}i\nu_0^2\\
i/\sqrt{3\alpha} & 2\nu_0\end{pmatrix}.
\end{align}
With this choice of $X$, it is clear that $\det(\frac{X}{z+\alpha}+I)=1$. It should be mentioned that we assume in \eqref{H} and \eqref{X} that $x$ does not belong to the zero set of the function $1-c(x)$ (see \eqref{c-simplify} below).

Since $\frac{X}{z+\alpha}+I$ is analytic at $z=0,2\alpha$, the parametrices near these points can be just defined by replacing $N(z)$ by $\widehat{N}(z)$ in \eqref{P0}, \eqref{E0} and \eqref{P+} and \eqref{E+}.

\subsection{Small norm problem}
With the modified model function $\widehat{N}(z)$ in \eqref{Nhat}, $\widehat{P}^{(-)}(z)$ in \eqref{P-m}, ${P}^{(0)}(z)$ in \eqref{P0} and ${P}^{(+)}(z)$ in \eqref{P+}, the final transformation is defined by
\begin{align}\label{D-to-Rhat}
\widehat R(z)=D(z)\left\{\begin{aligned}
&\widehat P^{(-)}(z)^{-1}, &&z\in U(-\alpha),\\
&P^{(0)}(z)^{-1}, &&z\in U(0),\\
&P^{(+)}(z)^{-1}, &&z\in U(2\alpha),\\
&\widehat N(z)^{-1}, &&\mathrm{elsewhere}.
\end{aligned}\right.
\end{align}
Therefore, we get the following RH problem for $\widehat R(z)$.
\begin{enumerate}[(a)]
\item $\widehat R(z)$ is analytic for $z\in \mathbb{C}\setminus\Gamma_R$, where $\Gamma_R$ is shown in Figure \ref{Gamma-R}.

\item $\widehat R(z)$ satisfies the jump relation $\widehat R_+(z)=\widehat R_-(z)J_{\widehat R}(z)$
\begin{equation}\label{Rhat-jump}
J_{\widehat R}(z)=\left\{\begin{aligned}
&\widehat P^{(-)}(z)\widehat N(z)^{-1}, &&z\in\partial U(-\alpha),\\
&P^{(0)}(z)\widehat N(z)^{-1}, &&z\in\partial U(0),\\
&P^{(+)}(z)\widehat N(z)^{-1}, &&z\in\partial U(2\alpha),\\
&\widehat N(z)J_D(z)\widehat N(z)^{-1}, &&z\in\Gamma_R\setminus(\partial U(-\alpha)\cup\partial U(0)\cup\partial U(2\alpha)),
\end{aligned}\right.
\end{equation}
where $J_D(z)$ denotes the jump matrix in \eqref{D-jump1} and \eqref{D-jump2}.

\item As $z\to\infty$, we have
  \begin{align}\label{Rhat-at-infty}
 \widehat R(z)=I+\frac{\widehat R_1}{z}+O(z^{-2}).
 \end{align}
\end{enumerate}

It is readily seen from \eqref{matching0}, \eqref{matching+} and \eqref{matching-m} that
\begin{align}\label{Rhatjump-estim}
J_{\widehat R}(z)=I+O(|x|^{-\frac54}), \qquad x\to-\infty.
\end{align}
Thus, we have that for $z\in\mathbb{C}\setminus\Gamma_R$,
\begin{align}\label{Rhat-estim}
\widehat R(z)=I+O(|x|^{-\frac54}), \qquad x\to-\infty.
\end{align}

\section{Proof of Theorem \ref{thm}}\label{Section5}
In this section, we prove  Theorem \ref{thm}. The proof depends on the formulas \eqref{y-express} and the asymptotic analysis of the RH problem for $\Psi(z,x)$ given in the previous two sections. Below, it should be understood that $|s_2|<1$ in Section \ref{proof1} and $|s_2|>1$ in Section \ref{proof2}.

\subsection{Derivation of \eqref{y-asymp1} and \eqref{H-asymp1}}\label{proof1}
By tracing back the transformations $\Psi(z,x)\to A(z)\to B(z)\to C(z)\to D(z)\to R(z)$ as given in \eqref{Psi-to-A}, \eqref{A-to-B}, \eqref{B-to-C}, \eqref{C-to-D} and \eqref{D-to-R}, we obtain that for $z\in\mathbb{C}\setminus\Gamma_R$,
\begin{align}\label{D=RN}
D(z)=R(z)N(z).
\end{align}
Taking the limit $z\to\infty$ in \eqref{D=RN}, it follows from the expansions \eqref{D-at-infty}, \eqref{N-expansion} and \eqref{R-at-infty} that
\begin{align}\label{B1-eq1}
B_{1,21}&=R_{1,21}+N_{1,21},\\
B_{1,22}&=R_{1,22}+N_{1,22}.\label{B1-eq2}
\end{align}
In views of \eqref{Psi1Psi2}-\eqref{H-express}, we obtain from \eqref{B1-21}, \eqref{B1-22}, \eqref{N1}, \eqref{B1-eq1} and \eqref{B1-eq2} that
\begin{align}\label{y-expre}
y(x)&=-\left(\alpha+4\sqrt{3\alpha}\,i\nu R_{1,21}+2R_{1,22}-R_{1,21}^2\right)|x|^{\frac12},\\
\mathcal{H}(x)&=-\frac{2\alpha}{3}|x|^{\frac32}
+2\sqrt{3\alpha}\,i\nu|x|^{\frac14}-R_{1,21}|x|^{\frac14}.\label{H-expre}
\end{align}

Next, we compute $R_1$. Since the RH problem for $R(z)$ is equivalent to the following singular integral equation
\begin{align}
R(z)=I+\frac{1}{2\pi i}\int_{\Gamma_R}\frac{R_-(\xi)(J_R(\xi)-I)d\xi}{\xi-z},
\end{align}
we get from \eqref{Rjump-estimate}, \eqref{R-estimate} and \eqref{PN-1} that as $x\to-\infty$,
\begin{align}
R_1&=-\frac{1}{2\pi i}\int_{\Gamma_R}R_-(\xi)(J_R(\xi)-I)d\xi\nonumber\\
&=-\frac{1}{2\pi i}\int_{\partial U(-\alpha)}(J_R(\xi)-I)d\xi+O(|x|^{-\frac{5}{4}})\nonumber\\
&=-\frac{|x|^{-\frac58}}{2\pi i}\int_{\partial U(-\alpha)}N(\xi)M(\xi)L(\xi)M(\xi)^{-1}N(\xi)^{-1}d\xi
+O(|x|^{-\frac{5}{4}}).
\end{align}
Via a direct residue calculation using \eqref{N}, \eqref{d}, \eqref{lambda-1-asymp} and \eqref{L}, we obtain that as $x\to-\infty$,
\begin{align}\label{R1-21}
R_{1,21}&=\frac{i}{2\sqrt{3\alpha}}(A-B)|x|^{-\frac58}+O(|x|^{-\frac{5}{4}}),\\
R_{1,22}&=\left[-\frac12(A+B)+\nu(A-B)\right]|x|^{-\frac58}+O(|x|^{-\frac{5}{4}}),
\label{R1-22}
\end{align}
where
\begin{align}\label{AB}
A=\frac{\nu s_2|x|^{\frac54\nu}e^{\frac{3\pi i\nu}{2}}2^{6\nu}(3\alpha)^{\frac{5\nu}{2}}}{h_1e^{2|x|^{\frac54}g_+(-\alpha)}
e^{-\frac{\pi i}4}2(3\alpha)^{\frac14}},\qquad
B=\frac{h_1e^{2|x|^{\frac54}g_+(-\alpha)}}{s_2|x|^{\frac54\nu}e^{\frac{3\pi i\nu}{2}}2^{6\nu}(3\alpha)^{\frac{5\nu}{2}}e^{-\frac{\pi i}4}2(3\alpha)^{\frac14}}.
\end{align}
Inserting the asymptotics in \eqref{R1-21}, \eqref{R1-22} into \eqref{y-expre} and \eqref{H-expre} gives
\begin{align}\label{y-expre1}
y(x)&=-\alpha|x|^{\frac12}+(A+B)|x|^{-\frac18}+O(|x|^{-\frac34}),\\
\mathcal{H}(x)&=-\frac{2\alpha}{3}|x|^{\frac32}
+2\sqrt{3\alpha}\,i\nu|x|^{\frac14}-\frac{i(A-B)}{2\sqrt{3\alpha}}|x|^{-\frac38}
+O(|x|^{-1}),\label{H-expre1}
\end{align}
as $x\to-\infty$.

Using \eqref{AB} and \eqref{h0h1}, we immediately have
\begin{align}
A\pm B=\frac{\nu\Gamma(-\nu) s_2|x|^{\frac54\nu}e^{\frac{\pi i\nu}{2}}2^{6\nu}(3\alpha)^{\frac{5\nu}{2}}}{\sqrt{2\pi}e^{2|x|^{\frac54}g_+(-\alpha)}
e^{-\frac{\pi i}4}2(3\alpha)^{\frac14}}\pm\frac{\sqrt{2\pi}h_1e^{2|x|^{\frac54}g_+(-\alpha)}}
{\Gamma(-\nu)s_2|x|^{\frac54\nu}e^{\frac{\pi i\nu}{2}}2^{6\nu}(3\alpha)^{\frac{5\nu}{2}}e^{-\frac{\pi i}4}2(3\alpha)^{\frac14}}.
\end{align}
With \eqref{Stokes-eq5}, \eqref{nu}, the relation $\Gamma(-\nu)=\overline{\Gamma(\nu)}$ and the reflection formula of the Gamma function
\begin{align}\label{Gamma-nu}
-\nu\Gamma(-\nu)\Gamma(\nu)=\frac{\pi}{\sin\pi\nu}=\frac{2\pi i}{e^{\pi i\nu}|s_2|^2},
\end{align}
we find that
\begin{align}\label{A+B}
A+B&=\frac{\sqrt{2\pi}}{(3\alpha)^{\frac14}|s_2||\Gamma(\nu)|e^{\frac{\pi i\nu}{2}}}\cos\varphi(x),\\
A-B&=\frac{i\sqrt{2\pi}}{(3\alpha)^{\frac14}|s_2||\Gamma(\nu)|e^{\frac{\pi i\nu}{2}}}\sin\varphi(x),
\end{align}
where
\begin{align}
\varphi(x)=2ig_+(-\alpha)|x|^{\frac54}-i\nu\ln\left[2^6(3\alpha)^{\frac52}|x|^{\frac54}
\right]-\arg\Gamma(\nu)+\arg s_2-\frac{\pi}{4}.
\end{align}
On the other hand, it follows from \eqref{Gamma-nu} that
\begin{align}\label{|Gamma-nu|}
|\Gamma(\nu)|=\frac{\sqrt{2\pi}}{\sqrt{i\nu}\,e^{\frac{\pi i\nu}{2}}|s_2|}.
\end{align}
Substituting \eqref{A+B}-\eqref{|Gamma-nu|} into \eqref{y-expre1} and \eqref{H-expre1} and then recalling \eqref{g}, \eqref{nu} and \eqref{g(-alpha)}, we arrive at \eqref{y-asymp1} and \eqref{H-asymp1} for the case $0<|s_2|<1$.

In the reduced case $|s_2|=0$, by applying the asymptotic analysis of the RH problem  for $\Psi$ discussed in Section \ref{sec:szero}, we obtain the following asymptotic formula
\begin{align}\label{y0}
y(x)&=-\left(-\frac{x}{6}\right)^{\frac12}+O(x^{-\frac34}),\\
\mathcal{H}(x)&=-4\left(-\frac{x}{6}\right)^{\frac32}+O(x^{-1}),\label{H0}
\end{align}
as $x\to-\infty$. We notice from \eqref{nu} that $\nu=0$ when $|s_2|=0$. Thus, \eqref{y0} and \eqref{H0} can be viewed respectively as the limits of \eqref{y-asymp1} and \eqref{H-asymp1} as $a\to0$.

\begin{rem}
Note that both the parabolic cylinder parametrix and Airy parametrix have a full asymptotic expansion as their variables tend to infinity. Hence, the jump matrices of $R(z)$ on $\partial U(-\alpha)$ and $\partial U(2\alpha)$ can also be fully asymptotically expanded as $|x|\to+\infty$. This indicates that we can theoretically obtain the full asymptotic expansion of the real solutions of Painlev\'{e} I equation. We conjecture that \eqref{y-asymp1} can be improved as
\begin{equation}\label{eq-full-expansion-PI}
y(x)\sim\sqrt{-\frac{x}{6}}\sum\limits_{k=0}^{\infty}\frac{a_{k}}{(-x)^{\frac{5k}{4}}}
+\frac{\sqrt{a}}{(-24x)^{\frac{1}{8}}}\sum\limits_{n=1}^{\infty}\cos(n\psi(x))
\sum\limits_{k=0}^{\infty}\frac{b_{k}^{(n)}}{(-x)^{\frac{5k}{8}}},
\end{equation}
where $\psi(x)=\frac{4\cdot24^{\frac14}}5(-x)^{\frac54}-\frac{5a}8\ln(-x)+\phi$ and $a_{k}, b_{k}^{(n)}$ are constants. One may try to determine the coefficients $a_{k}$ and $b_{k}^{(n)}$ by following the way above in this section, while we think it is better to obtain these coefficients by substituting \eqref{eq-full-expansion-PI} into the Painlev\'{e} I equation \eqref{PI} directly.
\end{rem}

\subsection{Derivation of \eqref{y-asymp2} and \eqref{H-asymp2}}\label{proof2}
Inverting the transformations $\Psi(z,x)\to A(z)\to B(z)\to C(z)\to D(z)\to\widehat{R}(z)$ as given in \eqref{Psi-to-A}, \eqref{A-to-B}, \eqref{B-to-C}, \eqref{C-to-D} and \eqref{D-to-Rhat}, we now get that for $z\in\mathbb{C}\setminus\Gamma_R$,
\begin{align}\label{D=RhatN}
D(z)=\widehat R(z)\widehat N(z).
\end{align}
Letting $z\to\infty$ in \eqref{D=RhatN}, we obtain from the expansions \eqref{D-at-infty}, \eqref{Nhat-expansion} and \eqref{Rhat-at-infty} that
\begin{align}\label{B1-eq1m}
B_{1,21}&=\widehat R_{1,21}+N_{1,21}+X_{21},\\
B_{1,22}&=\widehat R_{1,22}+N_{1,22}+X_{22}.\label{B1-eq2m}
\end{align}
Using \eqref{Psi1Psi2}-\eqref{H-express}, \eqref{B1-21}, \eqref{B1-22}, \eqref{N1}, \eqref{X}, \eqref{B1-eq1m} and \eqref{B1-eq2m}, we have
\begin{align}\label{y-expre-m}
y(x)&=-\alpha|x|^{\frac12}-\frac{12\alpha c}{(1-c)^2}|x|^{\frac12}-2\widehat R_{1,22}|x|^{\frac12}-\widehat R_{1,21}^2|x|^{\frac12}\nonumber\\
&\qquad\qquad\ -4\sqrt{3\alpha}\,i\nu_0 \widehat R_{1,21}|x|^{\frac12}
+2\sqrt{3\alpha}i\widehat R_{1,21}\frac{1-c}{1+c}|x|^{\frac12},\\
\mathcal{H}(x)&=-\frac{2\alpha}{3}|x|^{\frac32}
+2\sqrt{3\alpha}\,i\nu_0|x|^{\frac14}-\frac{i\sqrt{3\alpha}(1+c)}{1-c}|x|^{\frac14}
-\widehat R_{1,21}|x|^{\frac14}.\label{H-expre-m}
\end{align}

For the asymptotics of $\widehat R_1$, we use the singular integral equation for $\widehat R(z)$ and the estimations \eqref{Rhatjump-estim}, \eqref{Rhat-estim} to deduce that
as $x\to-\infty$,
\begin{align}\label{Rhat1-asym}
\widehat R_1=-\frac{1}{2\pi i}\int_{\Gamma_R}\widehat R_-(\xi)
(J_{\widehat R}(\xi)-I)d\xi=O(|x|^{-\frac{5}{4}}).
\end{align}
On the other hand, it follows from \eqref{c} and \eqref{h0h1} that
\begin{align}
c=\frac{\sqrt{2\pi}e^{-\frac12\pi i}e^{2|x|^{\frac54}g_+(-\alpha)}}
{\Gamma(\frac12-\nu_0)s_2|x|^{\frac54\nu_0}e^{\frac{\pi i\nu_0}{2}}2^{6\nu_0}(3\alpha)^{\frac{5\nu_0}{2}}}.
\end{align}
Since
\begin{align}
|\Gamma(1/2-\nu_0)|^2=\Gamma(1/2-\nu_0)\Gamma(1/2+\nu_0)
=\frac{\pi}{\sin\pi(\frac12-\nu_0)}=\frac{2\pi}{e^{\pi i\nu_0}|s_2|^2},
\end{align}
we have
\begin{align}\label{c-simplify}
c=\exp \left\{i\bigg(-2ig_+(-\alpha)|x|^{\frac54}+\frac54i\nu_0\ln|x|
\right.&+i\nu_0\ln[2^6(3\alpha)^{\frac52}]\nonumber\\
&\left.-\arg s_2-\arg\Gamma(\frac{1}{2}-\nu_0)-\frac\pi2\bigg)\right\}.
\end{align}
Hence, by \eqref{c-simplify},
\begin{align}\label{c-equa}
\frac{c}{(1-c)^2}=-\frac{1}{4\sin^2\omega(x)},\qquad\quad
\frac{1+c}{1-c}=-i\cot\omega(x),
\end{align}
where
\begin{align}\label{omega}
\omega(x)=ig_+(-\alpha)|x|^{\frac54}-\frac58i\nu_0\ln|x|&-
\frac12i\nu_0\ln[2^6(3\alpha)^{\frac52}]\nonumber\\
&+\frac12\arg s_2+\frac12\arg\Gamma(1/2-\nu_0)+\frac\pi4.
\end{align}
Plugging the asymptotics \eqref{Rhat1-asym}, the equations \eqref{c-equa} and \eqref{omega} into \eqref{y-expre-m}, \eqref{H-expre-m} and then recalling \eqref{g}, \eqref{nu} and \eqref{g(-alpha)}, we obtain \eqref{y-asymp2} and \eqref{H-asymp2}.

\section*{Acknowledgements}
The work of Jun Xia was supported in part by the National Natural Science Foundation of China [Grant no. 12401322] and the Guangdong Basic and Applied Basic Research Foundation (Grant no. 2024A1515012985).
The work of Wen-Gao Long was supported in part by the National Natural Science Foundation of China
[Grant no. 12401094], the Natural Science Foundation of Hunan Province [Grant no. 2024JJ5131] and the Outstanding Youth Fund of Hunan Provincial Department of Education [Grant no. 23B0454].

\begin{appendices}

\section{Model parametrices}\label{Appendix}
\subsection{Parabolic cylinder  parametrix}\label{PCP}
\begin{figure}[H]
  \centering
  \includegraphics[width=5cm,height=5cm]{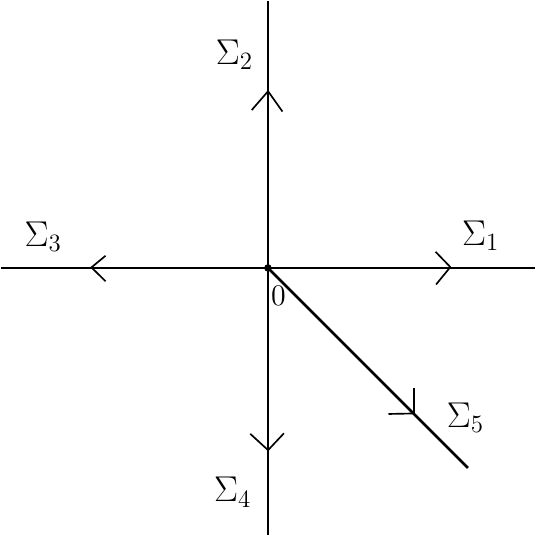}\\
  \caption{The jump contour $\Sigma_k$ }\label{PC}
\end{figure}

Let
$$\mathbf{P}(\lambda)=2^{-\frac{\sigma_{3}}{2}}\begin{pmatrix}
D_{-\nu-1}(i\lambda) & D_{\nu}(\lambda) \\D_{-\nu-1}'(i\lambda) & D_{\nu}'(\lambda)\end{pmatrix}
\begin{pmatrix}
e^{i \frac{\pi}{2}(\nu+1)} & 0 \\
0 & 1
\end{pmatrix},
$$
where $D_{\nu}(\lambda)$ is parabolic cylinder function (see \cite[Chapter 12]{NIST}).
Define
$$
H_{0}=\begin{pmatrix}1 & 0 \\ h_{0} & 1\end{pmatrix}, \
H_{1}=\begin{pmatrix}1 & h_{1} \\ 0 & 1\end{pmatrix}, \
H_{n+2}=e^{i \pi\left(\nu+\frac{1}{2}\right) \sigma_{3}} H_{n} e^{-i \pi\left(\nu+\frac{1}{2}\right) \sigma_{3}}, \ n=0,1,
$$
where
\begin{equation}\label{h0h1}
h_{0}=-i \frac{\sqrt{2 \pi}}{\Gamma(\nu+1)}, \quad h_{1}=\frac{\sqrt{2 \pi}}{\Gamma(-\nu)} e^{i \pi \nu}, \quad 1+h_{0} h_{1}=e^{2 \pi i \nu}.
\end{equation}
We introduce
$$
\Phi^{(\mathrm{PC})}(\lambda)=\left\{\begin{aligned}
&\mathbf{P}(\lambda),&& \arg\lambda \in\left(-{\pi}/{4}, 0\right), \\
&\mathbf{P}(\lambda)H_0,&& \arg\lambda \in\left(0,{\pi}/{2}\right), \\
&\mathbf{P}(\lambda)H_0H_1,&& \arg\lambda \in\left({\pi}/{2}, \pi\right), \\
&\mathbf{P}(\lambda)H_0H_1H_2,&& \arg\lambda \in\left(\pi,{3\pi}/{2}\right), \\
&\mathbf{P}(\lambda)H_0H_1H_2H_3,&& \arg\lambda \in\left({3\pi}/{2},{7\pi}/{4}\right).
\end{aligned}\right.
$$
Thus, $\Phi^{(\mathrm{PC})}$ solves the following RH problem (cf. \cite{FIKN}).
\begin{enumerate}[(a)]
\item $\Phi^\mathrm{(PC)}$ is analytic for $\lambda\in \mathbb{C}\setminus \bigcup^5_{k=1}\Sigma_{k}$, where $\Sigma_{k}=\{\lambda\in\mathbb{C}:\arg \lambda=(k-1)\pi/2\}$, $k=1,2,3,4$ and  $\Sigma_{5}=\{\lambda\in\mathbb{C}:\arg \lambda=-\pi/4\}$; see Figure \ref{PC}.
\item $\Phi^\mathrm{(PC)}$ satisfies the jump conditions
\begin{equation}\label{PC-jump}
\Phi^\mathrm{(PC)}_+(\lambda)=\Phi^\mathrm{(PC)}_-(\lambda)\left\{\begin{aligned}
&H_0, &&\lambda\in\Sigma_1,\\
&H_1, &&\lambda\in\Sigma_2,\\
&H_2, &&\lambda\in\Sigma_3,\\
&H_3, &&\lambda\in\Sigma_4,\\
&e^{2\pi i\nu\sigma_3}, && \lambda\in\Sigma_5.
\end{aligned}\right.
\end{equation}
\item As $\lambda\to\infty$, we have
\begin{align}\label{PC-at-infty}
\Phi^\mathrm{(PC)}(\lambda)=\begin{pmatrix}0 &1 \\ 1 & -\lambda\end{pmatrix}2^{\frac{\sigma_3}{2}}\left[I+\sum_{k=1}^{\infty}
\Phi^{\mathrm{(PC)}}_k\lambda^{-k}\right]
e^{\left(\frac{\lambda^{2}}{4}-\nu\ln \lambda\right) \sigma_{3}},
\end{align}
where $\Phi^{\mathrm{(PC)}}_{2k-1}$ is off-diagonal and $\Phi^{\mathrm{(PC)}}_{2k}$ is diagonal,
\begin{align}
\Phi^{\mathrm{(PC)}}_{1}=\begin{pmatrix}0&\nu \\1&0 \end{pmatrix},\qquad
\Phi^{\mathrm{(PC)}}_{2}=\begin{pmatrix}\frac{\nu(\nu+1)}{2}&0 \\0&-\frac{\nu(\nu-1)}{2} \end{pmatrix}.
\end{align}
\end{enumerate}

\subsection{Airy parametrix}\label{AP}
\begin{figure}[H]
  \centering
  \includegraphics[width=6cm,height=4cm]{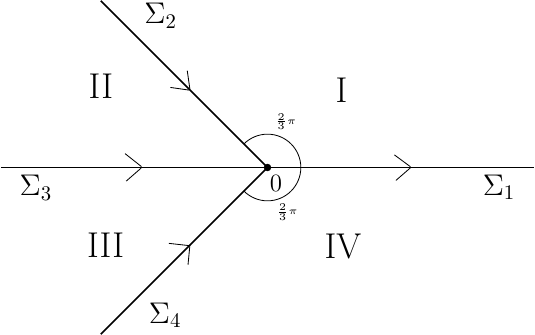}\\
  \caption{The jump contour $\Sigma_k$ and regions I,II,III,IV}\label{Airy}
\end{figure}

Let $\varphi=e^{2\pi i/3}$. We consider
\begin{equation}
\Phi^{(\mathrm{Ai})}(\lambda)=\Upsilon
\left\{\begin{aligned}
&\begin{pmatrix}
\mathrm{Ai}(\lambda) & \mathrm{Ai}(\varphi^{2}\lambda) \\
\mathrm{Ai}^{\prime}(\lambda) & \varphi^{2} \mathrm{Ai}^{\prime}(\varphi^{2}\lambda)
\end{pmatrix} e^{-i \frac{\pi}{6} \sigma_{3}}, &&\lambda\in \mathrm{I},\\
&\begin{pmatrix}
\mathrm{Ai}(\lambda) & \mathrm{Ai}(\varphi^{2}\lambda) \\
\mathrm{Ai}^{\prime}(\lambda) & \varphi^{2} \mathrm{Ai}^{\prime}(\varphi^{2}\lambda)
\end{pmatrix} e^{-i \frac{\pi}{6} \sigma_{3}}\begin{pmatrix}
1 & 0 \\
-1 & 1
\end{pmatrix},&&\lambda\in \mathrm{II},\\
&\begin{pmatrix}
\mathrm{Ai}(\lambda) & -\varphi^{2}\mathrm{Ai}(\varphi\lambda) \\
\mathrm{Ai}^{\prime}(\lambda) & -\mathrm{Ai}^{\prime}(\varphi\lambda)
\end{pmatrix} e^{-i \frac{\pi}{6} \sigma_{3}}\begin{pmatrix}
1 & 0 \\
1 & 1
\end{pmatrix},&&\lambda\in \mathrm{III},\\
&\begin{pmatrix}
\mathrm{Ai}(\lambda) & -\varphi^{2} \mathrm{Ai}(\varphi\lambda) \\
\mathrm{Ai}^{\prime}(\lambda) & -\mathrm{Ai}^{\prime}(\varphi\lambda)
\end{pmatrix} e^{-i \frac{\pi}{6} \sigma_{3}},&&\lambda\in \mathrm{IV},
\end{aligned}
\right.
\end{equation}
where $\mathrm{Ai}(\lambda)$ is the Airy function (cf. \cite[Chapter 9]{NIST}) and
$$
\Upsilon=\sqrt{2\pi}\begin{pmatrix}e^{\frac{1}{6}\pi i} & 0\\ 0 & e^{-\frac{1}{3}\pi i} \end{pmatrix}.
$$
It is direct to see that $\Phi^{(\mathrm{Ai})}$ satisfies the following RH problem (cf. \cite{Deift}).
\begin{enumerate}[(a)]
\item $\Phi^\mathrm{(Ai)}$ is analytic for $\lambda\in \mathbb{C}\setminus \bigcup^4_{k=1}\Sigma_{k}$, where $\Sigma_k$ is depicted in Figure \ref{Airy}.
\item $\Phi^\mathrm{(Ai)}$ satisfies the jump conditions
\begin{equation}\label{Airy-jump}
\Phi^\mathrm{(Ai)}_+(\lambda)=\Phi^\mathrm{(Ai)}_-(\lambda)\left\{
\begin{aligned}
&\begin{pmatrix}1 & 1 \\ 0 & 1 \end{pmatrix}, & \lambda\in\Sigma_{1},\\
&\begin{pmatrix}1 & 0 \\ 1 & 1 \end{pmatrix}, & \lambda\in\Sigma_{2}\cup\Sigma_{4},\\
&\begin{pmatrix}0 & 1 \\ -1 & 0 \end{pmatrix}, & \lambda\in\Sigma_{3}.
\end{aligned}\right.
\end{equation}
\item $\Phi^\mathrm{(Ai)}$ satisfies the asymptotic behavior
\begin{equation}\label{Airy-at-infty}
\Phi^{(\mathrm{Ai})}(\lambda)=\lambda^{-\frac{\sigma_{3}}{4} }\frac{1}{\sqrt{2}} \begin{pmatrix}1 & i\\ i &1\end{pmatrix}\left[I+O(\lambda^{\frac32})\right] e^{-\frac{2}{3}\lambda^{\frac{3}{2}} \sigma_{3}},\quad \lambda\to\infty.
\end{equation}
\end{enumerate}

\section{Stokes multipliers of a special solution of Painlev\'{e} I }\label{appendix2}
According to \cite{LongLiWang} (see also \cite{Kapaev-Kitaev-1993}), $\Psi(z,x)$ satisfies the following differential equation
\begin{equation}\label{eq-fold-Lax-pair}
\frac{\partial\Psi}{\partial z}
=\begin{pmatrix}y_{x}&2z^{2}+2yz+x+2y^2\\2(z-y)&-y_{x}\end{pmatrix}\Psi.
\end{equation}
The only singularity of the above equation is the irregular singular point at $z=\infty$.
Following \cite{Kapaev-Kitaev-1993} (see also \cite{Kap2}),
there exist the canonical solutions $\Psi_{k}(z,x)$, $k\in\mathbb{Z}$, of \eqref{eq-fold-Lax-pair} with the asymptotic expansion
\begin{equation}\label{eq-canonical-solutions}
\Psi^{(k)}(z,x)
=z^{\frac{1}{4}\sigma_{3}}\frac{\sigma_{3}+\sigma_{1}}{\sqrt{2}}
\left[I+\frac{-\mathcal{H}(x)\sigma_{3}}{\sqrt{z}}
+O\left(\frac{1}{z}\right)\right]
e^{(\frac{4}{5}\lambda^{\frac{5}{2}}+xz^{\frac{1}{2}})\sigma_{3}}
\end{equation}
as $z\rightarrow\infty$ with $z\in\Omega_{k}$, uniformly for all $x$ bounded away from the poles of the Painlev\'{e} I transcendents, where $\mathcal{H}(x)$ is introduced in \eqref{Hamilton}, and the canonical sectors are
$$\Omega_{k}=\left\{z\in\mathbb{C}:~\arg z\in \left(-\frac{3\pi}{5}+\frac{2k\pi}{5},\frac{\pi}{5}+\frac{2k\pi}{5}\right)\right\}, \qquad k\in\mathbb{Z}.$$
These canonical solutions are connected by
\begin{equation}\label{eq-Stokes-matrices}
\Psi^{(k+1)}(z,x)=\Psi^{(k)}(z,x)S_{k},\quad S_{2k-1}=\begin{pmatrix}1&s_{2k-1}\\0&1\end{pmatrix},\quad S_{2k}=\begin{pmatrix}1&0\\s_{2k}&1\end{pmatrix}.
\end{equation}
By \cite{Bertola-Tovbis-2013}, we can proceed further as follows.
Define
\begin{equation}\label{eq-def-hat{Phi}}
\hat{\Psi}(z,x)=\begin{pmatrix}0&1\\ 1&-\frac{1}{2}\left(-y_{x}+\frac{1}{2(z-y)}\right)\end{pmatrix}
(z-y)^{\frac{\sigma_{3}}{2}}\Psi(z,x)
\end{equation}
Then $\hat{\Phi}(z,x)$ satisfies
\begin{equation}\label{eq-system-hat-Phi}
\frac{\partial}{\partial z}\hat{\Psi}(z,x)=\begin{pmatrix}0&2\\V(z,x)&0\end{pmatrix}\hat{\Psi}(z,x),
\end{equation}
where
\[
2V(z,x)=y_{x}^{2}+4z^{3}+2z x-2y x-4y^{3}-\frac{y_{x}}{z-y}+\frac{3}{4}\frac{1}{(z-y)^2}.
\]
Moreover, using the expansion of $\Psi(z,x)$ in \eqref{eq-canonical-solutions}, one can verify that
\begin{equation}\label{eq-canonical-solutions-hat-Phi}
\hat{\Psi}^{(k)}(z,x)
=\frac{z^{-\frac{3}{4}\sigma_{3}}}{\sqrt{2}}
\begin{pmatrix}1&-1\\1&1\end{pmatrix}
\left(I+O(z^{-\frac{1}{2}})\right)
e^{(\frac{4}{5}z^{\frac{5}{2}}+xz^{\frac{1}{2}})\sigma_{3}}
\end{equation}
as $z\rightarrow\infty$ with $z\in\Omega_{k}$.
The asymptotic expansions of $\Psi^{(k)}(z,x)$ and $\hat{\Psi}^{(k)}(z,x)$ in \eqref{eq-canonical-solutions} and \eqref{eq-canonical-solutions-hat-Phi} are valid
only when $x$ is bounded away from $p$.
When $x\to p$, according to \cite{Bertola-Tovbis-2013},
the system \eqref{eq-system-hat-Phi} reduces to
\begin{equation}\label{eq-system-hat-Phi2}
\frac{\partial}{\partial z}\hat{\Psi}(z,p)
=\begin{pmatrix}0&2\\2z^3+pz-14H&0\end{pmatrix}\hat{\Psi}(z,p),
\end{equation}
and the corresponding asymptotic expansions of $\hat{\Psi}_{k}(z,p)$ in \eqref{eq-canonical-solutions-hat-Phi} should be replaced by
\begin{equation}\label{eq-canonical-solutions-near-pole}
\hat{\Psi}^{(k)}(z,p)
=\frac{z^{-\frac{3}{4}\sigma_{3}}}{\sqrt{2}}
\begin{pmatrix}-i&-i\\-i&i\end{pmatrix}
\left(I+O\left(z^{-\frac12}\right)\right)
e^{(\frac{4}{5}z^{\frac{5}{2}}+pz^{\frac{1}{2}})\sigma_{3}}
\end{equation}
as $z\rightarrow\infty$ with $z\in\Omega_{k}$; see \cite[Corollary A.8]{Bertola-Tovbis-2013}. It is readily seen that $\hat{\Psi}^{(k)}$ also satisfy the connection formula
\begin{equation}\label{eq-connection-hat-Psi}
\hat{\Psi}^{(k+1)}(z,p)=\hat{\Psi}^{(k)}(z,p)S_{k}, k\in\mathbb{Z}.
\end{equation}
Finally, if denoting
\[
\hat{\Psi}(z,p)=\begin{pmatrix}\phi_{1}\\\phi_{2}\end{pmatrix},
\]
and letting $Y(z;p)=\phi_{1}$, then we arrive at the following reduced triconfluent Heun equation (RTHE); see  \cite[Eq. (6)]{Xia-Xu-Zhao}
\begin{equation}\label{Schrodinger-equation-triconfluent-Heun}
\frac{d^{2}Y}{d z^{2}}=\left[4z^{3}+2pz-28H\right]Y.
\end{equation}

When $(p,H)=(0,0)$, the equation \eqref{Schrodinger-equation-triconfluent-Heun} can be explicitly solved by the Bessel functions. For any solution of \eqref{Schrodinger-equation-triconfluent-Heun}, there exist two constants $C_{1}$ and $C_{1}$ such that
\begin{equation}
Y(z)=C_{1}z^{\frac{1}{2}}I_{\frac{1}{5}}\left(\frac{4}{5}z^{\frac{5}{2}}\right)
+C_{2}z^{\frac{1}{2}}K_{\frac{1}{5}}\left(\frac{4}{5}z^{\frac{5}{2}}\right).
\end{equation}
Applying the asymptotic behaivor of the Bessel functions $I_{\nu}(z)$ and $K_{\nu}(z)$ (cf. see \cite[Eqs. 10.40.1, 10.40.2]{NIST}) and matching them with the behaviors of $(\hat{\Psi}^{(k)})_{11}$ and $(\hat{\Psi}^{(k)})_{12}$ in \eqref{eq-canonical-solutions-near-pole}, we have
\begin{equation}\label{eq-Bessel-asym--pi/5}
\begin{split}
\left(z^{\frac{1}{2}}I_{\frac{1}{5}}\left(\frac{4}{5}z^{\frac{5}{2}}\right),
z^{\frac{1}{2}}K_{\frac{1}{5}}\left(\frac{4}{5}z^{\frac{5}{2}}\right)\right)&\sim ((\Psi^{(0)})_{11},(\Psi^{(0)})_{12})
\left(\begin{matrix}
c_{1}&0\\
c_{2}&c_{3}
\end{matrix} \right)
\end{split}
\end{equation}
as $z\to\infty$ with $\arg{z}\sim -\frac{\pi}{5}$, where
\begin{equation}\label{eq-def-c1-c2-c3}
c_{1}=\left(\frac{4\pi}{5}\right)^{-\frac{1}{2}}i,\quad c_{2}=\left(\frac{4\pi}{5}\right)^{-\frac{1}{2}}e^{-\frac{\pi i}{5}},\quad c_{3}=\left(\frac{4}{5\pi}\right)^{-\frac{1}{2}}i.
\end{equation}
Similarly when $z\to\infty$ with $\arg{z}\sim \frac{\pi}{5}$, we have
\begin{equation}\label{eq-Bessel-asym-pi/5}
\begin{split}
\left(z^{\frac{1}{2}}I_{\frac{1}{5}}\left(\frac{4}{5}z^{\frac{5}{2}}\right),
z^{\frac{1}{2}}K_{\frac{1}{5}}\left(\frac{4}{5}z^{\frac{5}{2}}\right)\right)&\sim ((\Psi^{(1)})_{11},(\Psi^{(1)})_{12})
\left(\begin{matrix}
c_{1}&0\\
\tilde{c}_{2}&c_{3}
\end{matrix} \right),
\end{split}
\end{equation}
where $\tilde{c}_{2}=-e^{\frac{2\pi i}{5}}c_{2}=ie^{\frac{\pi i}{5}}c_{1}$.
Hence, according to the connection formula \eqref{eq-connection-hat-Psi}, we conclude that
\begin{equation}
S_{0}=\left(\begin{matrix}
c_{1}&0\\
c_{2}&c_{3}
\end{matrix} \right)\left(\begin{matrix}
c_{1}&0\\
\tilde{c}_{2}&c_{3}
\end{matrix} \right)^{-1}
=\left(\begin{matrix}
1&0\\
\frac{c_{2}-\tilde{c}_{2}}{c_{1}}&1
\end{matrix} \right),
\end{equation}
which implies that $s_{0}=\frac{c_{2}-\tilde{c}_{2}}{c_{1}}=-2i\cos\left(\frac{\pi}{5}\right)$.

Combining \cite[Eqs. 10.40.1, 10.40.2]{NIST} with the connection formulas between $I_{\nu}(z)$ and $K_{\nu}(z)$ (cf. see \cite[Eqs. 10.34.1, 10.34.2]{NIST}), we can further obtain
\begin{equation}\label{eq-Bessel-asym-3pi/5}
\begin{split}
\left(z^{\frac{1}{2}}I_{\frac{1}{5}}\left(\frac{4}{5}z^{\frac{5}{2}}\right),
z^{\frac{1}{2}}K_{\frac{1}{5}}\left(\frac{4}{5}z^{\frac{5}{2}}\right)\right)
\sim\left((\Psi^{(2)})_{11},(\Psi^{(2)})_{12}\right)\left(\begin{matrix}
d_{1}&d_{3}\\d_{2}&d_{4}
\end{matrix}\right),
\end{split}
\end{equation}
as $z\to\infty$ with $\arg{z}\sim \frac{3\pi}{5}$,
where
\begin{equation}
\label{eq-def-d1-d2-d3-d4}
d_{1}=ie^{\frac{\pi i}{5}}\tilde{c}_{2}=-c_{1},\quad  d_{2}=ie^{\frac{\pi i}{5}}c_{1}=-c_{2},\quad d_{3}=e^{-\frac{\pi i}{5}}ic_{3}+\pi \tilde{c}_{2},\quad d_{4}=\pi c_{1}=c_{3}.
\end{equation}
Comparing \eqref{eq-Bessel-asym-pi/5} and \eqref{eq-Bessel-asym-3pi/5}, it is readily seen that
\begin{equation}
S_{1}=\left(\begin{matrix}
c_{1}&0\\
\tilde{c}_{2}&c_{3}
\end{matrix} \right)\left(\begin{matrix}
d_{1}&d_{3}\\d_{2}&d_{4}
\end{matrix}\right)^{-1}=\left(\begin{matrix}
1&-c_{1}d_{3}\\0&1
\end{matrix}\right).
\end{equation}
This implies that $s_{1}=-c_{1}d_{3}=-2i\cos{\left(\frac{\pi}{5}\right)}$.

Finally, making use of \eqref{Stokes-eq1}, we find that $s_{k}=-2i\cos\left(\frac{\pi}{5}\right)$ for all $k\in\mathbb{Z}$.

\end{appendices}


\end{document}